\documentclass[12pt]{article}   %12pt --> art10

\usepackage{amsfonts}

\newcommand{\cM}{{\cal M}}
\newcommand{\cP}{{\cal P}}
\newcommand{\cF}{{\cal F}}
\newcommand{\cO}{{\cal O}}
\newcommand{\cN}{{\cal N}}
\newcommand{\wJ }{{\cal J}}

\newcommand{\Sym}{{\rm Sym}}
\newcommand{\Stab}{{\rm Stab}}
\newcommand{\rank}{{\rm rank}}
\newcommand{\coker}{{\rm coker}}
\newcommand{\Max}{{\rm Max}}
\newcommand{\Gal}{{\rm Gal}}
\newcommand{\ord}{{\rm ord}}
\renewcommand{\char}{{\rm char}}

\renewcommand{\ker}{{\rm ker}}
\newcommand{\End}{{\rm End}}
\newcommand{\Hom}{{\rm Hom}}
\newcommand{\id}{{\rm id}}
\newcommand{\Spec}{{\rm Spec}}

\newcommand{\ZZ}{{\mathbb Z}}
\newcommand{\CC}{{\mathbb C}}
\newcommand{\RR}{{\mathbb R}}
\newcommand{\NN}{{\mathbb N}}
\renewcommand{\AA}{{\mathbb A}}

\newcommand{\QQ}{{\mathbb Q}}
\newcommand{\HH}{{\mathbb H}}

\newcommand{\wq}{{\mathfrak q}}
\newcommand{\wP}{{\mathfrak P}}
\newcommand{\wA}{{\mathfrak A}}
\renewcommand{\wp}{{\mathfrak p}}
\renewcommand{\wJ}{{\mathfrak J}}
\newcommand{\wD}{{\mathfrak D}}

\newcommand{\Gammam}{{\mit\Gamma}}

\newcommand{\fis}{{\#}}
\newcommand{\RG}{R\fis\Gammam}
\newcommand{\KnullRG}{K_0(\RG)}
\newcommand{\ra}{\rightarrow}
\newcommand{\okg}{\cO_K\Gammam}
\newcommand{\okpg}{\cO_{K_\wp}\Gammam}
\newcommand{\Kpq}{{\bar{K}_\wp}}
\newcommand{\Kpqg}{{\Kpq\Gammam}}
\newcommand{\MtM}{\overline{\cM_i \otimes \cM_i}}
\newcommand{\bigwedgem}{\Lambda}
\def\rightepi{{\longrightarrow \kern-0.7em \rightarrow}}
\newcommand{\oplusm}{\mathop{\oplus}\limits}

\begin{document}

\vspace*{15ex}

\begin{center}
{\LARGE\bf Operations on Locally Free Classgroups}\\
\bigskip
by\\
\bigskip
{\sc Bernhard K\"ock}
\end{center}

\bigskip

\begin{quote}
{\footnotesize {\bf Abstract.} 
Let $\Gammam$ be a finite group and $K$ a number field.
We show that the operation $\psi_k$ defined by Cassou-Nogu\`es and
Taylor on the locally free classgroup ${\rm Cl}(\cO_K\Gammam)$ is
a symmetric power operation if $\gcd(k,\ord(\Gammam)) = 1$. Using the
equivariant Adams-Riemann-Roch theorem, we furthermore give a geometric
interpretation of a formula established by Burns and Chinburg for 
these operations.}

\end{quote}

\bigskip

\section*{Introduction}

Exterior or symmetric powers of modules 
play an important role in algebraic $K$-theory,
in representation theory of finite groups, in number theory, or in 
algebraic geometry. The object of this paper is to work out their central
role in a fascinating interplay of these four subjects we are now going
to describe.

The $k$-th Adams operation $\psi^k$ on the classical ring $K_0(\CC\Gammam)$
of virtual characters of a finite group $\Gammam$ can be defined in the
following two ways: In terms of modules, $\psi^k$ is given by 
\[\psi^k(M) = N_k(\bigwedgem^1_\CC(M), \ldots, \bigwedgem^k_\CC(M));\]
i.e., the exterior powers $\bigwedgem^1_\CC(M), \ldots,
\bigwedgem^k_\CC(M)$ of the representation $M$ of $\Gammam$ are plugged into
the $k$-th Newton polynomial $N_k = t_1^k + \cdots + t_k^k$ considered
as a polynomial in the elementary symmetric functions of $t_1, \ldots, t_k$. 
In terms of characters, $\psi^k$ is given by 
\[\psi^k(\chi)(\gamma) = \chi(\gamma^k)\]
(for $\gamma \in \Gammam$ and $\chi \in K_0(\CC\Gammam)$ a character).

Now let $K$ be a number field. Let $K_0(\cO_K\Gammam)$ denote the 
Grothendieck group of all projective modules over the group ring 
$\cO_K\Gammam$ and let
\[{\rm Cl}(\cO_K\Gammam) = \ker(\rank: K_0(\cO_K\Gammam) \ra \ZZ)\]
denote the locally free classgroup associated with $K$ and $\Gammam$. 
Fr\"ohlich has established an isomorphism between ${\rm Cl}(\cO_K\Gammam)$
and a certain residue class group of the group
\[\Hom_{{\rm Gal}(\bar{K}/K)}(K_0(\bar{K}\Gammam), {\cal I}(\bar{K}))\]
consisting of
Galois-invariant homomorphisms from $K_0(\bar{K}\Gammam) = 
K_0(\CC\Gammam)$ to the group 
${\cal I}(\bar{K})$ of 
ideles of $\bar{K}$ (e.g., see Theorem 1 on p.\ 20 in 
\cite{F1}). Using Taylor's group
logarithm techniques, Cassou-Nogu\`es and Taylor have shown that, under
certain additional assumptions, the endomorphism $\Hom(\psi^k, {\cal I}
(\bar{K}))$ of $\Hom_{{\rm Gal}(\bar{K}/K)} 
(K_0(\bar{K}\Gammam),{\cal I}(\bar{K}))$
induces an endomorphism on ${\rm Cl}(\cO_K\Gammam)$ (e.g., see
Chapter 9 in \cite{T}). This endomorphism will be denoted by
$\psi_k^{\rm CNT}$ in the following considerations. In the paper \cite{BC},
Burns and Chinburg have studied the question whether there is an algebraic
description of $\psi_k^{\rm CNT}$, for instance in terms of power 
operations on $\cO_K\Gammam$-modules. Their main result essentially is
a formula for $\psi_k^{\rm CNT}(\wA)$ where $\wA$ is an arbitrary 
$\Gammam$-stable
fractional ideal in a tame Galois extension $N/K$ with Galois group
$\Gammam$. Their formula expresses $\psi_k^{\rm CNT}(\wA)$ in terms of
powers of the different $\wD_{N/K}$ and of powers of the ideal $\wA$
(see Corollary 2.4 in \cite{BC}).

The central result of this paper is the identification of the endomorphism
$\psi_k^{\rm CNT}$ with a symmetric power operation on $K_0(\cO_K\Gammam)$
in the case $\gcd(k,\ord(\Gammam))=1$. The aim of this paper finally is 
to give a geometric interpretation of the formula of Burns and Chinburg,
namely as an equivariant Adams-Riemann-Roch formula for the $\Gammam$-morphism
$\Spec(\cO_N) \ra \Spec(\cO_K)$. 

The following fact (see Proposition 1.1) 
is the fundamental observation for our central result:
Let $R$ be a commutative ring and $M$ a projective module over the group
ring $R\Gammam$; if $\gcd(k, \ord(\Gammam))$ is invertible in $R$, then
the $k$-th exterior power $\bigwedgem^k_R(M)$ and the $k$-th 
symmetric power $\Sym^k_R(M)$ are projective over $R\Gammam$ again; here,
as usual, $\Gammam$ acts on the exterior or symmetric power diagonally.

A slightly more general form of this observation will enable us to define a 
$k$-th symmetric power operation $\sigma^k$ on the Grothendieck group
$K_0(R\Gammam)$ and, more generally, on Quillen's higher 
$K$-groups $K_q(R\Gammam)$, $q \ge 0$, (see section 1). An obvious generalization
of the module-theoretic definition of $\psi^k$ on $K_0(\CC\Gammam)$
explained above furthermore yields an (additive) $k$-th Adams operation
$\psi^k$ on $K_q(R\Gammam)$, $q\ge 0$, if $\gcd(k, \ord(\Gammam))$ is
invertible in $R$ (see section 1). In Theorem 1.6, we will prove several
standard properties of
these operations such as e.g.\ the compatibility with the Cartan homomorphism.
Moreover, Theorem 1.6 contains a proof of the identity
\[\psi^k([R\Gammam]) = [R\Gammam]\]
conjectured in (1.12) in Chapter 9 of \cite{T}.

In section 2, we will show that $\psi^k$ is multiplicative. Since there
is no splitting principle available, this cannot be proved as e.g.\ in
\cite{FL}. Instead, we will use the universal form of the 
Cauchy decomposition and of the Pieri formula (established by Akin,
Buchsbaum, and Weyman in characteristic-free representation theory) in 
order to express the exterior power $\bigwedgem^i(M\otimes N)$ of the
tensor product of two modules as a universal polynomial in the 
exterior powers of $M$ and $N$. It remains an open problem how to
prove that $\psi^k \circ \psi^l = \psi^{kl}$ on $K_0(R\Gammam)$ (here,
$l$ denotes another natural number such that $\gcd(l, \ord(\Gammam))$
is invertible in $R$).

We now assume that $\gcd(k, \ord(\Gammam)) = 1$ and fix a $k' \in \NN$
with $kk' \equiv 1$ mod $\ord(\Gammam)$. By the constructions explained above,
we have operations $\sigma^k$ and $\psi^k$ on $K_0(\cO_K\Gammam)$ and,
in particular, on ${\rm Cl}(\cO_K\Gammam)$ by restricting. The precise
form of our central result now is:
\[\psi_{k'}^{\rm CNT} = \sigma^k \quad {\rm and} \quad 
k \cdot \psi_{k'}^{\rm CNT} = \psi^k \quad {\rm on} \quad 
{\rm Cl}(\cO_K\Gammam)\]
(see Theorem 3.7). This result in particular answers the question
(raised in (1.12) in Chapter 9 of \cite{T}) how
the endomorphism $\psi_k^{\rm CNT}$ behaves with respect to the 
Cartan homomorphism (see Corollary 3.9). 

Now let $N/K$ be a tame Galois extension with Galois group $\Gammam$. In
section 4, we will essentially establish the following simple representation
of the cotangential element $T_f^\vee$ associated with the 
$\Gammam$-morphism $f: \Spec(\cO_N) \ra \Spec(\cO_K)$:
\[T_f^\vee = [\wD_{N/K}^{-1}] - [\cO_N].\]
Here, as above, $\wD_{N/K}$ denotes the different associated with $N/K$.

This representation of $T_f^\vee$ implies a particularly simple 
shape of the Bott multiplier and then of the equivariant 
Adams-Riemann-Roch formula for $f$ (see Theorem 5.4). On the other hand,
we will reformulate the formula of Burns and Chinburg mentioned above using
the central result Theorem 3.7 (see Theorem 5.6). It will then turn out 
that, more or less, the formula of Burns and Chinburg is a strengthening
of the equivariant Adams-Riemann-Roch formula.

{\bf Acknowledgments.} First of all, I would like to thank C.~Greither
for inviting me to stay with him at Universit\'e Laval in Qu\'ebec for
a few weeks in September/October 1996. In discussions during this visit,
he made the observation Proposition 1.1 which has later turned out to
be the most important tool in understanding the operations on locally 
free classgroups defined by Cassou-Nogu\`es and Taylor. Moreover, he has 
contributed several algebraic arguments which have considerably improved the
paper. It is a pleasure for me to thank him for all this assistance. Furthermore,
I would like to thank J.~Weyman for the reference to the universal form
of the Pieri formula and for long e-mail discussions about the 
plethysm problem. Finally, I would like to thank D.~Burns who asked
me some time ago how to geometrically interpret the formula he developed
jointly with T.~Chinburg. 

\bigskip

{\bf Notations.} Let $\Gammam$ be a finite group and let $R$ be a commutative
$\Gammam$-ring, i.e., a commutative ring together with an action of $\Gammam$
by ring automorphisms. The twisted group ring associated with $\Gammam$ and
$R$ is denoted by $\RG$. If $\Gammam$ acts on $R$ trivially, we simply 
write $R\Gammam$ for $\RG$. We denote the Grothendieck group of all
finitely generated projective $\RG$-modules (of all finitely generated 
$R$-projective $\RG$-modules, of all finitely generated $\RG$-modules) by
$\KnullRG$, (by $K_0(\Gammam, R)$, and by $K'_0(\RG)$, respectively). We have
a canonical homomorphism $c: \KnullRG \ra K_0(\Gammam, R)$ which is called
the Cartan homomorphism. If the order of the group $\Gammam$ is invertible
in $R$, $c$ is bijective (e.g., see the proof of Corollary 2.2c) in 
\cite{KoAdHi}) and $\KnullRG$ will be identified with $K_0(\Gammam, R)$. 
Furthermore, we have a canonical homomorphism $K_0(\Gammam, R) \ra
K'_0(\RG)$. If $R$ is regular (e.g., if $R$ is a field or a ring of 
integers in a number field), this homomorphism is bijective (e.g., see
Satz (2.1) in \cite{KoARR}) and $K_0(\Gammam,R)$ will be identified with
$K'_0(\RG)$. Thus, if $R$ is a field of characteristic zero, all three
Grothendieck groups coincide. The corresponding higher $K$-groups
will be denoted by $K_q(\RG)$, $K_q(\Gammam, R)$, $K'_q(\RG)$, $q\ge 0$, 
(see \cite{Q}). The corresponding identifications hold for higher
$K$-groups as well. For $K_1(\RG)$, we will use also the following two descriptions
(e.g., see \cite{B}): By the third theorem on p.\ 228 in \cite{Gr1}, 
$K_1(\RG)$ is isomorphic to the abelianized infinite linear group 
${\rm GL}_\infty (\RG)^{\rm ab}$; furthermore, $K_1(\RG)$ is the
free Abelian group over all pairs $(M,\alpha)$, where $M$ is
a f.~g.\ projective $\RG$-module and $\alpha$ is an $\RG$-automorphism
of $M$, modulo the relations defined on p.\ 348 in \cite{B}. The notation
$K_0T(-)$ will always stand for a Grothendieck group of torsion modules.

\bigskip

\section*{\S 1 Power Operations for Projective Modules over \\
Group Rings}

Let $\Gammam$ be a finite group, $R$ a commutative $\Gammam$-ring and 
$k$ a natural number such that ${\rm gcd}(k, \ord(\Gammam))$ is invertible
in $R$.

The aim of this section is to construct a symmetric power operation 
$\sigma^k$, an exterior power operation $\lambda^k$, and an Adams 
operation $\psi^k$ on Quillen's higher $K$-groups $K_q(\RG)$, $q\ge 0$,
in particular on the Grothendieck group $\KnullRG$ of all f.\ g.\
projective $\RG$-modules. Furthermore, we will prove several standard
properties of these operations and we will show the equality
$\psi^k([\RG]) = [\RG]$ if $\gcd(k,\ord(\Gammam)) = 1$. The following 
proposition is the fundamental 
observation for this equality and for the whole section. 

{\bf Proposition 1.1.} Let $k_1, \ldots, k_r \in \NN$ such that
$\gcd(k_1, \ldots, k_r, \ord(\Gammam))$ is invertible in $R$. (For
instance, this holds if $k_1 + \ldots + k_r = k$ or if at least one 
$k_i$ is $1$.) Let $M_1, \ldots, M_r$ be projective 
$\RG$-modules. Then the $\RG$-module 
\[ P:= \Sym^{k_1}(M_1) \otimes \cdots \otimes \Sym^{k_r}(M_r)\]
is projective over $\RG$, too. (Here, all tensor products and all 
symmetric powers in $P$ are formed over $R$ and $\Gammam$ acts on
$P$ diagonally.) \\
If moreover $\gcd(k_1, \ldots, k_r, \ord(\Gammam)) = 1$ and $M_1, \ldots, M_r$ are
free over $\RG$, then $P$ is free over $\RG$, too.

{\bf Proof.} First, let $r=1$. We write $k$ for $k_1$. Let $M$ 
be a free $\RG$-module and $B \subset M$ an $R$-basis of $M$ which
carries a free action of $\Gammam$. For instance, the basis consisting of
the $\Gammam$-images of an $\RG$-basis of $M$ is such a basis. Let 
$B^k/\Sigma_k$ denote the set of unordered $k$-tuples in $B$. For
any $H = (b_1, \ldots, b_k) \in B^k/\Sigma_k$, let $b_H := \prod_{j=1}^k
b_j$ denote the corresponding standard basis element of $\Sym^k(M)$. 
Then we obviously have $\gamma(b_H) = b_{\gamma(H)}$ for all $\gamma \in
\Gammam$, where the action of $\Gammam$ on $B^k/\Sigma_k$ is defined in
the obvious way. For any $H\in B^k/\Sigma_k$, the stabilizer
$\Stab(H):= \{\gamma \in \Gammam: \gamma(H) = H\}$ clearly acts on
the sets $H(n):= \{b \in B: b \mbox{ occurs } n \mbox{ times in } H\}$,
$n\ge 1$, and this action is free by assumption. 
Thus, the order $\ord(\Stab(H))$ divides the number 
$\fis  H(n)$ of elements in $H(n)$ for all $n$. Hence, it divides 
$k = \sum_{n=1}^\infty n\cdot (\fis  H(n))$, too.\\
Now, let $r$ be arbitrary and 
$M_1, \ldots, M_r$ free $\RG$-modules. Furthermore, let $B_1, \ldots,
B_r$ be bases of $M_1, \ldots, M_r$ as above. Then, for any 
$H = (H^1, \ldots, H^r) \in B_1^{k_1}/\Sigma_{k_1} \times \cdots \times
B_r^{k_r}/\Sigma_{k_r}$, the element $b_H := b_{H^1} \otimes \cdots \otimes
b_{H^r}$ is a standard basis element of $P$ with $\gamma(b_H) = 
b_{\gamma(H)}$. The order $\ord(\Stab(H))$ of the
stabilizer of $H$ divides the order of the stabilizer of $H^i$ for all $i$.
Then, by the first part of the proof, it divides 
$\gcd(k_1, \ldots, k_r, \ord(\Gammam))$, too. 
Hence, $\ord(\Stab(H))$ is invertible in $R$. Thus, the canonical 
$\RG$-linear surjection $\RG \ra R\fis [\Gammam/\Stab(H)]$, $\gamma \mapsto
[\gamma \Stab(H)]$, can be split by the well-defined $\RG$-linear map
\[ R\fis [\Gammam/\Stab(H)] \ra \RG, \quad [\gamma \Stab(H)] \mapsto 
\frac{1}{\ord(\Stab(H))} \sum_{\gamma' \in \Stab(H)} [\gamma \gamma'].\]
This proves that the $\RG$-module $R\fis [\Gammam/\Stab(H)]$ is projective
over $\RG$. The $\RG$-submodule generated by $b_H$ obviously is isomorphic
to $R\fis [\Gammam/\Stab(H)]$. Thus, it is $\RG$-projective, too. Being a
direct sum of such modules, also $P$ is projective. This shows the first
assertion in Proposition 1.1 for free $\RG$-modules $M_1, \ldots, M_r$
and, by passing to direct summands, also for projective $\RG$-modules.\\
If finally $\gcd(k_1,\ldots, k_r, \ord(\Gammam))=1$, then $\Stab(H)$ is trivial. 
Hence, $P$ is a free $\RG$-module.

Now we are going to construct the power operations mentioned above. For 
this, it is convenient to define the following exact subcategories
$\cM_i$, $i=1,2, \ldots$ of the exact category $\cM$ of all f.\ g.\
$R$-projective $\RG$-modules: Let $\cM_i$ be the smallest full subcategory
of $\cM$ which is closed under extensions and kernels of $\RG$-epimorphisms
and which contains all modules of the form $\Sym^{k_1}(M_1) \otimes \cdots
\otimes \Sym^{k_r}(M_r)$, where $M_1, \ldots, M_r$ are f.~g.\ projective
$\RG$-modules and $k_1, \ldots, k_r$ are natural numbers with
$k_1 + \cdots + k_r = i$. We call the modules in $\cM_i$ {\em modules
of weight $i$}. Obviously, $\cM_1$ is the category of all f.~g.\ 
projective $\RG$-modules. By Proposition 1.1, the category $\cM_k$ is
contained in $\cM_1$. 

{\bf Lemma 1.2.} For all $i, j \in \NN$, the tensor product induces
a functor
\[\otimes : \cM_i \times \cM_j \ra \cM_{i+j}.\]

{\bf Proof.} A module of weight $i$ is called of level $0$ iff it is of
the form $\Sym^{k_1}(M_1) \otimes \cdots \otimes \Sym^{k_r}(M_r)$ as above.
Inductively, a module of weight $i$ is called of level $\le m$ iff it
is an extension or the kernel of an $\RG$-epimorphism of two modules
of weight $i$ and level $\le m-1$. Then one easily shows by induction
that the tensor product of a module of weight $i$ and level $ \le m$ with
a module of weight $j$ and level $\le n$ is a module of weight $i+j$ 
and level $\le m+n$. Since any module of weight $i$ or $j$ is of finite
level, this shows Lemma 1.2. 

For any exact category $\cP$ and for any $i \ge 1$, let $\cF_i(\cP)$ denote
the category of all sequences $M_1 \hookrightarrow \cdots \hookrightarrow
M_i$ of admissible monomorphisms in $\cP$ of length $i$. Here, a 
monomorphism $M' \ra M$ in $\cP$ is called {\em admissible} iff it can be 
completed to a short exact sequence $0 \ra M' \ra M \ra M'' \ra 0$ in 
$\cP$. 

{\bf Lemma 1.3.} For any sequence $M_1 \hookrightarrow \cdots 
\hookrightarrow M_i$ in $\cF_i(\cM_1)$, the $\RG$-module 
\[M_1   \cdots   M_i := {\rm Image}(M_1 \otimes \cdots \otimes
M_i \, \, \stackrel{\rm can}{\longrightarrow} \, \, \Sym^i(M_i))\]
is contained in $\cM_i$.

{\bf Proof.} In the canonical filtration 
\[ M_1   \cdots   M_i \subseteq M_1  \cdots   M_{i-2} 
M_{i}  M_i \subseteq \cdots \subseteq \Sym^i (M_i)\]
of $\Sym^i(M_i)$, each successive quotient is isomorphic to one of the
modules $M_1   \cdots   M_{r-1} \otimes \Sym^{i-r+1}\left(
\frac{M_i}{M_r}\right)$, $r=1, \ldots, i$. Now, Lemma 1.3 follows from
Lemma 1.2 by induction on $i$. 

It is easy to prove that the exact categories $\cM_1, \cM_2, \ldots$ 
together with the tensor products $\otimes : \cM_i \times \cM_j \ra \cM_{i+j}$, 
$i,j \in \NN$, (well-defined by Lemma 1.2) and the symmetric power
operations $\cF_i(\cM_1) \ra \cM_i$, $(M_1 \hookrightarrow \cdots
\hookrightarrow M_i) \mapsto M_1   \cdots   M_i$, (well-defined
by Lemma 1.3) satisfy the axioms (E1) through (E5) in section 7 in \cite{Gr2}.
Unlike Grayson in section 7 in \cite{Gr2}, we cannot define operations
$\cF_i(\cM_j) \ra \cM_{ij}$ for $j > 1$, since it is not clear whether
the $i$-th symmetric power $\Sym^i(M)$ of a module $M$ of weight $j$ is a
module of weight $ij$. This problem is related to the so-called 
plethysm problem (see also Remark 2.6). Nevertheless, we may
apply the construction of section 7 in \cite{Gr2} for $j=1$. It yields
continuous maps 
\[\sigma^i: |G\cM_1| \ra |G^i\cM_i|, \quad i \ge 1,\]
from the geometric realization of the $G$-construction associated with 
$\cM_1$ to the geometric realization of the $i$-fold iterated 
$G$-construction associated with $\cM_i$. By passing to homotopy groups,
we obtain {\em symmetric power operations}
\[\sigma^i: K_q(\RG) \ra K_q(\cM_i), \quad q \ge 0, \quad i\ge 1,\]
on Quillen's $K$-groups. For $i=k$, we finally obtain a symmetric power operation
\[\sigma^k: K_q(\RG) \ra K_q(\RG)\]
on the higher $K$-groups $K_q(\RG)$, $q\ge 0$, since $\cM_k$ is
contained in $\cM_1$ by Proposition 1.1.

{\bf Lemma 1.4.} For any sequence $M_1 \hookrightarrow \cdots
\hookrightarrow M_i$ in $\cF_i(\cM_1)$, the $\RG$-module 
\[M_1 \wedge \ldots \wedge M_i := {\rm Image}(M_1 \otimes \cdots \otimes
M_i \, \, \stackrel{\rm can}{\longrightarrow} \, \, \bigwedgem^i(M_i))\]
is contained in $\cM_i$.

{\bf Proof.} For any $i \ge 1$ and any f.~g.\ projective $\RG$-module
$M$, the Koszul complex
\[0 \ra \bigwedgem^i(M) \ra M\otimes \bigwedgem^{i-1}(M) \ra \cdots \ra
\Sym^{i-1}(M) \otimes M \ra \Sym^i(M) \ra 0\]
is an exact sequence of $\RG$-modules. It follows from Lemma 1.2 by 
induction on $i$ that $\bigwedgem^i(M)$ is contained in $\cM_i$. Now,
Lemma 1.4 can be deduced from this fact as in the proof of Lemma 1.3.

{\bf Remark 1.5.} Lemma 1.4, Lemma 1.2, and Proposition 1.1 imply that 
Proposition 1.1
holds as well if some symmetric powers are replaced by exterior powers
in the definition of $P$. This can of course also be shown directly as in 
Proposition 1.1 by introducing some signs. Moreover, it is more convenient
to define the categories $\cM_i$, $i \ge 1$, with symmetric powers than
with exterior powers since it would otherwise not be clear whether 
symmetric powers (and, more generally, (co)Schur modules, see Proposition 2.1) 
are contained in the categories $\cM_i$, $i\ge 1$.

As above, we now obtain {\em exterior power operations}
\[ \lambda^i : |G\cM_1| \ra |G^i \cM_i|, \quad i \ge 1,\]
and the corresponding maps on the $K$-groups. 

By Lemma 1.2, we furthermore have continuous maps 
\[\otimes: |G^i \cM_i| \times |G^j\cM_j| \ra |G^{i+j} \cM_{i+j}|, 
\quad i,j \ge 1,\]
(see also section 1 of \cite{KoSh}). Let $N_i(s_1,\ldots,s_i)$ denote the $i$-th
Newton polynomial $t_1^i + \cdots + t_i^i$ considered as a polynomial
in the elementary symmetric functions $s_1, \ldots, s_i$ of $t_1, \ldots, 
t_i$. Replacing
$s_1, \ldots, s_i$ by the exterior power operations $\lambda^1, \ldots, 
\lambda^i$ defined above and the products in $N_i(s_1, \ldots, s_i)$ 
by the tensor products defined above, we finally obtain continuous maps
\[\psi^i : |G\cM_1| \ra |G^i\cM_i|, \quad i \ge 1,\]
and the corresponding maps on the $K$-groups as above. These maps are
called {\em Adams operations}.

Applying Grayson's construction to the exact category $\cM$ consisting of
all f.~g.\ $R$-projective $\RG$-modules (i.e., all categories 
$\cM_i$, $i\ge 1$, in section 7 in \cite{Gr2} are identical to $\cM$),
we similarly obtain power operations
\[\sigma^i, \lambda^i, \psi^i: K_q(\Gammam, R) \ra K_q(\Gammam, R), \quad
i \ge 1, \quad q \ge 0,\]
(see also \S 2 and \S 3 in \cite{KoARR}). 

{\bf Theorem 1.6.} \\
(a) Let $\gamma$ be one of the operations $\sigma^i$, $\lambda^i$, 
$\psi^i$, $i\ge 1$. Then the following diagram commutes for all $q\ge 0$:
\[\begin{array}{ccc}
K_q(\RG) & \stackrel{c}{\longrightarrow} & K_q(\Gammam, R)\\
\\
\downarrow \gamma && \downarrow \gamma \\
\\
K_q(\cM_i) & \stackrel{\rm can}{\longrightarrow} & K_q(\Gammam, R).
\end{array}\]
In particular, the operations $\sigma^k$, $\lambda^k$, and $\psi^k$ on 
$K_q(\RG)$ are compatible with the corresponding operations on 
$K_q(\Gammam, R)$ with respect to the Cartan homomorphism (for all 
$q \ge 0$). \\
(b) The operations $\sigma^i$, $\lambda^i$, $\psi^i$, $i \ge 1$, 
commute with base change with respect to any homomorphism $R\ra R'$
of $\Gammam$-rings and they commute with the restriction map for any
subgroup $\Gammam'$ of $\Gammam$.\\
(c) For all $i \ge 1$ and for all $M, N \in \cM_1$ we have in $K_0(\cM_i)$:\\
\hspace*{1em}(i) $\sigma^i([M] - [N]) = 
\sum_{{a \ge 0, \, b_1, \ldots, b_u \ge 1}\atop{a+b_1+\cdots+b_u = i}} 
(-1)^u [\Sym^a(M)\otimes \Sym^{b_1}(N) \otimes \cdots
\otimes \Sym^{b_u}(N)]$.\\
\hspace*{1em}(ii) The assertion (i) analogously holds for exterior powers.\\
\hspace*{1em}(iii) $\psi^i([M] + [N]) = \psi^i([M]) + \psi^i([N])$.\\
\hspace*{1em}(iv) $\sum _{j=0}^i (-1)^j \lambda^j([M]) 
\cdot \sigma^{i-j}([M]) = 0$.\\
(d) For all $i \ge 1$ and $q \ge 1$, we have:\\
\hspace*{1em}(i) The maps $\lambda^i, \sigma^i, \psi^i : K_q(\RG) \ra K_q(\cM_i)$
are homomorphisms.\\
\hspace*{1em}(ii) $\sigma^i = (-1)^{i-1} \lambda^i : K_q(\RG) \ra K_q(\cM_i)$.\\
\hspace*{1em}(iii) $\psi^i = (-1)^{i-1} i \lambda^i :  K_q(\RG) \ra K_q(\cM_i)$.\\
(e) If $\gcd(k, \ord(\Gammam)) = 1$, we have:
\[\psi^k([\RG]) = [\RG] \quad {\rm in} \quad \KnullRG.\]
(f) Let $\tilde{K}_0(\ZZ, \RG) := \ker(K_0(\ZZ, \RG) \,\, 
\stackrel{\rm can}{\longrightarrow} \, \, K_0(\RG))$ denote the
reduced Grothen\-dieck group of all pairs $(M,\alpha)$, where $M$ is
a f.~g.\ projective $\RG$-module and $\alpha$ is an $\RG$-automorphism
of $M$. Let $\gamma: \tilde{K}_0(\ZZ, \RG) \ra \tilde{K}_0(\ZZ, \RG)$
be one of the operations $\sigma^k$, $\lambda^k$, $\psi^k$ defined
similarly as above. Then the following diagram commutes:
\[\begin{array}{ccc}
\tilde{K}_0(\ZZ,\RG) & \stackrel{\rm can}{\rightepi} & K_1(\RG)\\
\\
\downarrow \gamma && \downarrow \gamma \\
\\
\tilde{K}_0(\ZZ,\RG) & \stackrel{\rm can}{\rightepi} & K_1(\RG).
\end{array}\]

{\bf Proof.} The assertions (a) and (b) follow from the functoriality of
Grayson's construction. \\
The assertions (c)(i) and (c)(ii) are proved in section 8 in \cite{Gr2}.
The $\lambda$-series $\lambda_t([M]) := \sum_{i \ge 0} \lambda^i([M]) t^i$
satisfies the rule $\lambda_t([M+N]) = \lambda_t([M]) \cdot \lambda_t([N])$
by (c)(ii); being the logarithmic derivative of the $\lambda$-series (see
page 23 in \cite{FL}), the $\psi$-series is additive. This proves 
assertion (c)(iii). The assertion (c)(iv) follows from the fact that 
the Koszul complex
\[0 \ra \bigwedgem^i(M) \ra M\otimes \bigwedgem^{i-1}(M) \ra
\cdots \ra \Sym^{i-1}(M) \otimes M \ra \Sym^i(M) \ra 0\]
is an exact sequence of $\RG$-modules. \\
The assertion (d)(i) follows from the fact that maps on higher homotopy
groups which are induced by continuous maps are homomorphisms. By the 
same argument, the products $K_q(\cM_i) \times K_q(\cM_j) \ra
K_q(\cM_{i+j})$ induced by the tensor product $\cM_i \times \cM_j \ra \cM_{i+j}$ 
are linear and hence 
trivial for $q \ge 1$. This together with Newton's recursion formula 
proves assertion (d)(iii).
Assertion (d)(ii) can be proved as in Theorem 6.1 in \cite{KoSh}.\\ 
To show assertion (e), we express $\psi^k$ as a polynomial in the
symmetric power operations $\sigma^1, \ldots, \sigma^k$ using (c)(iv).
Then all monomials
in this polynomial are of weight $k$. Thus, by Proposition 1.1, there
is an $n\in \NN$ such that $\psi^k([\RG]) = n \cdot [\RG]$. We have $n=1$
since $\rank_R(\psi^k([\RG])) = \rank_R([\RG])$. This proves assertion (e).\\
The assertion (f) can be proved as in section 9 in \cite{Gr2} (see also the
proof of Theorem 3.3 in \cite{KoAdHi}). 

{\bf Remark 1.7.} \\
(a) If not only $\gcd(k, \ord(\Gammam))$, but 
$k$ is invertible in $R$, we have constructed an Adams operation
$\psi^k$ on $K_q(\RG)$, $q\ge 0$, in \cite{KoAdPro} by using so-called
cyclic powers. One should be able to prove that that operation agrees
with the Adams operation defined in this section. Since there is no 
splitting principle available for $K_0(\RG)$, this is a problem already
(and in fact mainly) for $K_0(\RG)$. At the end of the paper \cite{KoAdPro},
some speculations are given from which a solution for this problem would
follow. If $R$ is a field and $\Gammam$ acts on $R$ trivially, the
Cartan homomorphism $c : K_0(R\Gammam) \ra K_0(\Gammam, R)$ is injective 
(see \S 5 in \cite{Ke}). Since both operations $\psi^k$ commute with $c$, they
agree in this case. In particular, Proposition 1.1 yields a considerable
simplification of \S 5 in \cite{Ke}. \\
(b) Compared with the construction of this section, the construction in
\cite{KoAdPro} has the disadvantage that it is only possible under
the stronger assumption $k\in R^\times$. However, it has the advantage that both the
multiplicativity of $\psi^k$ and the rule $\psi^k \circ \psi^l = \psi^l \circ
\psi^k$ can rather easily be proved and that the definition of $\psi^k([M])$
involves only two modules if $k$ is a prime number.

{\bf Remark 1.8.} Let $\Gammam'$ be a  subgroup of $\Gammam$. 
In \S 6 of \cite{KoGRR}, the
following question has been raised: Under which conditions does the
equality $\psi^k([R\fis  [\Gammam/\Gammam']]) = [R\fis [\Gammam/\Gammam']]$ hold
in $K_0(\Gammam, R)$. Using the equivariant Adams-Riemann-Roch formula,
we have proved  in \S 6 in \cite{KoGRR} that this equality holds in
a certain completion of $K_0(\Gammam, R)[k^{-1}]$ for all $k\in \NN$. 
By Theorem 1.6(e), we
now have that this equality holds already in $K_0(\Gammam, R)$ if $\Gammam' $ is
a normal subgroup and $\gcd(k, \ord(\Gammam/\Gammam'))=1$.

\bigskip

\section*{\S 2 Multiplicativity of Adams Operations}

Let $\Gammam$ be a finite group, $R$ a commutative $\Gammam$-ring, and
$k\in \NN$ such that $\gcd(k,\ord(\Gammam))$ is invertible in $R$.

The aim of this section is to prove the multiplicativity of the Adams
operation $\psi^k$ on $K_0(\RG)$ (and on the higher $K$-groups) defined
in the previous section. Since there is no splitting principle available
for $\KnullRG$, this cannot be proved as e.g.\ in \cite{FL}. In place
of the splitting principle, we will use a filtration of the exterior
power $\bigwedgem^i(M \otimes N)$ of the tensor product of two modules
$M$ and $N$ which has been constructed by K.\ Akin, D.\ Buchsbaum,
and J.\ Weyman. Each successive quotient in this filtration is
isomorphic to a tensor product of (co)Schur modules (see \cite{ABW}). 
Furthermore, we will use a universal form of the Pieri formula 
established by K.\ Akin and D.\ Buchsbaum in \cite{AB} to express
(co)Schur modules as polynomials in exterior powers. The results 
of this section will not be used in the subsequent sections.

Let $\lambda = (\lambda_1 \ge \lambda_2 \ge \ldots) \in \NN^\infty$ be
a partition with $1 \le |\lambda| := \lambda_1 + \lambda_2 + \cdots 
< \infty$. The transposed partition (see p.\ 217 in \cite{ABW}) is
denoted by $\tilde{\lambda} = (\tilde{\lambda}_1 \ge \tilde{\lambda}_2 \ge
\ldots)$. For any f.~g.\ free $R$-module $F$, the Schur module
$L_\lambda(F)$ is defined as the image of a certain canonical map
\[ d_\lambda(F): \bigwedgem^{\lambda_1}(F) \otimes \bigwedgem^{\lambda_2}(F)
\otimes \cdots \ra \Sym^{\tilde{\lambda}_1}(F) \otimes 
\Sym^{\tilde{\lambda}_2}(F) \otimes \cdots \]
(see pages 219 and 220 in \cite{ABW}). By Theorem II.2.16 in \cite{ABW},
$L_\lambda(F)$ is a f.~g.\ free $R$-module and $L_\lambda$ commutes
with base change. For any f.~g.\ projective $R$-module $M$, we define the
Schur module $L_\lambda(M)$ in the same way. By localization,
we see that $L_\lambda(M)$ is a f.~g.\ projective $R$-module. If 
furthermore $M$ is an $\RG$-module, then the map $d_\lambda(M)$ 
obviously is $\RG$-linear and $L_\lambda(M)$ is an $\RG$-module. The
same holds for the coSchur module $K_\lambda(M)$ defined on p.\ 220
in \cite{ABW}. 

Let $s_\lambda$ denote the Schur polynomial associated with $\tilde{\lambda}$
(see I \S 3 in \cite{MacD}). It will be viewed as a polynomial in the
elementary symmetric functions. It is homogeneous of weight $|\lambda|$.
Thus, for any f.~g.\ projective $\RG$-module $M$, 
\[s_\lambda(M) := s_\lambda([\bigwedgem^1(M)], [\bigwedgem^2(M)], \ldots)\]
is a well-defined element of $K_0(\cM_{|\lambda|})$ by Lemma 1.2 and
Lemma 1.4. (Here, the category $\cM_{|\lambda|}$ is the category
of modules of weight $|\lambda|$ defined in section 1.)

{\bf Proposition 2.1.} Let $M$ be a f.~g.\ projective $\RG$-module.
Then the Schur module $L_\lambda(M)$ and the coSchur module $K_\lambda(M)$
are contained in the category $\cM_{|\lambda|}$. Furthermore, we have:
\[[L_\lambda(M)] = s_\lambda(M) \quad {\rm and} \quad [K_\lambda(M)] =
s_{\tilde{\lambda}}(M) \quad {\rm in} \quad K_0(\cM_{|\lambda|}).\]

{\bf Proof.} First, we show that the Schur module $L_\lambda(M)$ is
contained in $\cM_{|\lambda|}$. For this, we use induction on $\lambda$
with respect to lexicographic order (i.e., $\lambda < \mu$ iff
there is an $r\in \NN$ such that $\lambda_1 = \mu_1, \ldots,
\lambda_r = \mu_r$ and $\lambda_{r+1} < \mu_{r+1}$). For $\lambda
= (1, 0, 0, \ldots)$, the assertion is trivial. For the induction step,
let $\mu$ denote that partition which arises from $\lambda$ by 
omitting the last
column and let $p$ denote the size of the last column in $\lambda$
(i.e., $p:= \tilde{\lambda}_{l(\tilde{\lambda})}$ and 
$\mu_i := \lambda_i -1$ for $i=1, \ldots, p$ and $\mu_i := \lambda_i$
for $i >p$). By the universal form of the Pieri formula (see 
Theorem (3)(a) in section 3 in \cite{AB}), we have
for all f.~g.\ free $R$-modules $F$: There is a filtration
on the tensor product $L_\mu(F) \otimes \Sym^p(F)$ such that the
successive quotients are isomorphic to Schur modules. More precisely,
we have: If $\nu$ runs through the set
\[I:=\{\nu \mbox{ a partition}: \mu_i \le \nu_i \le \mu_i + 1 \mbox{ for all } i
\mbox{ and } |\nu| = |\lambda|\},\]
then $L_\nu(F)$ runs through the sequence of successive quotients. The
partition $\lambda$ is obviously the biggest partition in $I$ with 
respect to lexicographic order; in particular, $L_\lambda(F)$ is
isomorphic to the smallest (nontrivial) module in the filtration of
$L_\mu(F) \otimes \Sym^p(F)$ (see also section 2 in \cite{AB}). By
localization, one easily sees that the corresponding assertions hold
for f.~g.\ $R$-projective $\RG$-modules as well. By the induction
hypothesis and by Lemma 1.2, the module 
$L_\mu(M)\otimes \Sym^p(M)$ and the modules
$L_\nu(M)$, $\nu \in I \backslash \{\lambda\}$, are of weight $|\lambda|$.
Thus, all filtration steps and, in particular, $L_\lambda(M)$ are
of weight $|\lambda|$. Using the corresponding theorem for coSchur
modules, one similarly deduces from this that $K_{|\lambda|}(M)$ is
of weight $|\lambda|$. \\
In order to show the equality $[L_\lambda(M)] = s_\lambda(M)$, we replace
the Pieri formula for $L_\mu(M) \otimes \Sym^p(M)$ by
the Pieri formula for $L_\mu(M) \otimes \bigwedgem^p(M)$ (see Theorem (3)(b)
in \cite{AB}) and the lexicographic order by the transposed lexicographic
order in the procedure above. Then, we immediately obtain a polynomial
$s_\lambda'$ such that $[L_\lambda(M)] = 
s_\lambda'([\bigwedgem^1(M)], [\bigwedgem^2(M)], \ldots)$
in $K_0(\cM_{|\lambda|})$. We have $s_\lambda' = s_\lambda$ by the Pieri
formula for Schur polynomials (cf. Appendix A in \cite{FH}).
Thus, we have shown the desired equality $[L_\lambda(M)] =
s_\lambda(M)$ in $K_0(\cM_{|\lambda|})$. The equality $[K_\lambda(M)] =
s_{\tilde{\lambda}}(M)$ can be shown analogously.

For any $i \ge 1$, let $\MtM$ be the smallest subcategory of $\cM$
which is closed under extensions and kernels of $\RG$-epimorphisms
and which contains all modules of the form $M\otimes N$, where $M,N \in \cM_i$.
The category $\MtM$ is obviously contained in $\cM_{2i}$. By
Proposition 1.1, $\overline{\cM_1 \otimes \cM_1}$ is contained in
$\cM_1$. Thus, by Proposition 1.1, $\overline{\cM_k \otimes \cM_k}$ is
contained in $\cM_1$. In particular, the tensor product defines 
a ring structure on $\KnullRG$ (in general without $1$). Let $P_i$ be
the universal polynomial defined on page 5 in \cite{FL}. It is
homogeneous of weight $i$ in both sets of variables. Thus, for all
f.~g.\ projective $\RG$-modules $M$, $N$, 
\[P_i(M,N) := P_i([\bigwedgem^1(M)], \ldots, [\bigwedgem^i(M)]; 
[\bigwedgem^1(N)], \ldots, [\bigwedgem^i(N)])\]
is a well-defined element of $K_0(\MtM)$ by Lemma 1.2 and Lemma 1.4. 

{\bf Theorem 2.2.} For all f.~g.\ projective $\RG$-modules $M$, $N$ and
all $i\ge 1$, the exterior power $\bigwedgem^i(M\otimes N)$ is 
contained in $\MtM$ and we have:
\[[\bigwedgem^i(M\otimes N)] = P_i(M,N) \quad {\rm in} \quad K_0(\MtM).\]

{\bf Proof.} By Theorem III.2.4 in \cite{ABW}, there is a filtration
on $\bigwedgem^i(M\otimes N)$ by $\RG$-submodules with the following
property: If $\lambda$ runs through the set of partitions of weight $i$,
then the tensor product $L_\lambda(M) \otimes K_\lambda(M)$ runs through
the sequence of successive quotients. (This is proved in \cite{ABW} only
in the non-equivariant situation and only for free $R$-modules $M$, $N$.
The general assertion follows from this using the standard arguments 
already explained above.) Hence, by Proposition 2.1, all filtration
steps and, in particular, $\bigwedgem^i(M\otimes N)$, are contained in
$\MtM$ and we have:
\begin{eqnarray*}
\lefteqn{[\bigwedgem^i(M\otimes N)] = \sum_{|\lambda| = i} [L_\lambda(M) 
\otimes K_\lambda(N)]}\\
&&= \sum_{|\lambda|= i} s_\lambda(M) \cdot s_{\tilde{\lambda}}(N)\\
&&= P_i(M,N) \quad {\rm in} \quad K_0(\MtM).
\end{eqnarray*}
Here, the last equality is formula (4.3') on p.\ 35 in \cite{MacD}.

{\bf Corollary 2.3.} For all $x, y \in \KnullRG$ and all $i \ge 1$,
we have
\[ \psi^i(x\cdot y) = \psi^i(x) \cdot \psi^i(y) \quad {\rm in} \quad
K_0(\MtM).\]
In particular, $\psi^k : \KnullRG \ra \KnullRG$ is a ring homomorphism.

{\bf Proof.} It obviously suffices to show Corollary 2.3 for elements
$x= [M]$ and $y=[N]$ where $M$, $N$ are f.~g.\ projective $\RG$-modules.
In this case, Corollary 2.3 follows from Theorem 2.2 and the following
identity of polynomials which is easy to prove:
\[N_i(P_j(X_1, \ldots, X_j; Y_1, \ldots, Y_j), j=1, \ldots, i) =
N_i(X_1, \ldots, X_i) \cdot N_i(Y_1, \ldots, Y_i).\]

The tensor product yields a $\KnullRG$-module structure on Quillen's
$K$-groups $K_q(\RG)$, $q\ge 0$, in the usual way (e.g., see \S 3 in 
\cite{Q}). 

{\bf Corollary 2.4.} Let $q \ge 0$. For all $x \in \KnullRG$ and all
$y \in K_q(\RG)$, we have:
\[\psi^k(x \cdot y) = \psi^k(x) \cdot \psi^k(y) \quad {\rm in} \quad
K_q(\RG).\]

{\bf Proof.} For $q=0$, this has been proved in Corollary 2.3. So let
$q \ge 1$. For any arbitrary group $G$ let $K_0(G, \RG)$ denote the
Grothendieck group of all representations of the group $G$ on f.\ g.\
projective $\RG$-modules. As in section 1, one can define a $k$-th
Adams operation $\psi^k$ on $K_0(G, \RG)$. Using Quillen's universal
transformation,
we then obtain a $k$-th Adams operation $\psi^k$ on $K_q(\RG)$ 
as in section 2 of \cite{KoAdHi}. As in section 9 of \cite{Gr2}, one
easily shows that this Adams operation agrees with the Adams operation
defined in section 1. Furthermore, one easily proves as in Corollary 2.3
that $\psi^k$ on $K_0(G, \RG)$ is multiplicative. Finally, as in
Corollary 2.2a) in \cite{KoAdHi}, one deduces Corollary 2.4 from this
using Quillen's universal transformation.

{\bf Remark 2.5.} If $\Gammam$ acts on $R$ trivially, the tensor product
defines a $K_0(\Gammam, R)$-module structure on $K_0(R\Gammam)$ by Lemma 2.1
in \cite{KoAdPro}. One easily shows as above that the $k$-th Adams operation
$\psi^k: K_0(R\Gammam) \ra K_0(R\Gammam)$ is semi-linear with respect to
the $k$-th Adams operation on $K_0(\Gammam, R)$ defined e.g.\ in section 1.

{\bf Remark 2.6.} Let $l$ be another natural number such that 
$\gcd(l, \ord(\Gammam))$ is invertible in $R$. One should be able to prove
similarly to Corollary 2.3 that $\psi^k \circ \psi^l = \psi^{kl}$ in 
$\End(\KnullRG)$. For this, one should be able to construct a ``universal
decomposition'' of $\bigwedgem^k (\bigwedgem^l(M))$ similarly to the
universal form of the
Cauchy decomposition of $\bigwedgem^i(M\otimes N)$ used in the proof of
Theorem 2.2. This problem is called the {\em plethysm problem} in the
literature and is not yet solved.\\
If $R$ is a field and $\Gammam$ acts on $R$ trivially, then the
Cartan homomorphism $c: K_0(R\Gammam) \ra K_0(\Gammam,R)$ is injective (see
\S 5 in \cite{Ke}). Thus, the rule $\psi^k \circ \psi^l = \psi^{kl}$ for
$K_0(R\Gammam)$ follows from the corresponding rule for $K_0(\Gammam, R)$ 
in this case. If
$R$ is the ring of integers in a number field and $\Gammam$ acts on
$R$ trivially, the rule $\psi^k \circ \psi^l = \psi^{kl}$ follows from 
Theorem 1.3 (1) on p.\ 99 in \cite{T} and Theorem 3.7 which is the main  
result of the next section. 

\bigskip

\section*{\S 3 Power Operations on Locally Free Classgroups}

Let $K$ be a number field, $\Gammam$ a finite group, and $k\in \NN$ such that
$k$ is coprime to the exponent $e(\Gammam)$ of $\Gammam$.
Furthermore, we fix a $k' \in \NN$ with
$k k' \equiv 1 $ mod $ e(\Gammam)$. Let $\cO_K$ denote the ring 
of integers in $K$ and ${\rm Cl}(\okg):= \ker(\rank: K_0(\okg) \ra \ZZ)$
the locally free classgroup associated with $K$ and $\Gammam$.

In this section, we will prove that the endomorphism $\psi_{k'}^{\rm CNT}$ of
${\rm Cl}(\okg)$ defined by Ph.\ Cassou-Nogu\`es and M.\ Taylor in
\cite{CNT} coincides with the restriction of the symmetric power operation 
$\sigma^k$ on $K_0(\okg)$ defined in section 1.

To begin with, we recall some basic facts about the locally free
classgroup (e.g., see \S 1 in Chapter 1 in \cite{T}). Let 
$K_0T(\okg)$ denote the Grothendieck group of all $\cO_K$-torsion modules
over $\okg$ which have a $\okg$-projective resolution of finite length
(and then also of length $1$ by Theorem 8 on p.\ 145 in \cite{Se2}). Let
$\partial: K_1(K\Gammam) \ra K_0T(\okg)$ be the connecting homomorphism;
it is uniquely determined by the following property: Let $F$ be a f.\ g.\
free $\okg$-module and let $\alpha \in {\rm GL}(F\otimes K) \cap 
\End(F)$; then $\partial((F\otimes K, \alpha)) = \coker(\alpha: F\ra F)$.
Furthermore, let $i_*: K_0T(\okg) \ra K_0(\okg)$ denote the resolution map:
it maps a torsion module $M$ with the $\okg$-projective resolution
$0 \ra P \ra Q \ra M \ra 0$ to the element $[P] - [Q]$ in $K_0(\okg)$. Then
the sequence
\[ K_1(K\Gammam) \,\, \stackrel{\partial}{\longrightarrow} \,\,
K_0T(\okg) \,\, \stackrel{i_*}{\longrightarrow} \,\,
K_0(\okg) \,\, \stackrel{\rank}{\longrightarrow} \,\, \ZZ\]
is exact; in particular, we have:
\[{\rm Cl}(\okg) \cong \coker(\partial: K_1(K\Gammam) \ra 
K_0T(\okg)).\]
Decomposing torsion modules into $\wp$-primary torsion modules induces
an isomorphism
\[K_0T(\okg) \cong \oplusm_{\wp \in \Max(\cO_K)} K_0T(\okpg);\]
here, $K_\wp$ denotes the completion of $K$ at $\wp$ and
$K_0T(\okpg)$ is the Grothendieck group of all torsion modules over 
$\okpg$ which have a $\okpg$-free resolution of length $1$. Finally,
we have canonical exact (localization) sequences
\[K_1(\okpg) \longrightarrow K_1(K_\wp\Gammam) \,\,
\stackrel{\partial}{\longrightarrow} \,\,
K_0T(\okpg) \longrightarrow 0, \quad \wp \in \Max(\cO_K);\]
here, $K_1(\okpg) \ra K_1(K_\wp \Gammam)$ is the canonical map and the local
connecting homomorphism $\partial: K_1(K_\wp\Gammam) \ra K_0T(\okpg)$
is defined in the same way as the global one. 

Now we fix a prime ideal $\wp \not= 0$ in $\cO_K$. Let $\Phi$ denote
the composition
\[\Phi: K_1(K_\wp\Gammam) \,\, 
\stackrel{\partial}{\longrightarrow} \,\,
K_0T(\okpg) \,\, \stackrel{\rm can}{\longrightarrow} \,\,
K_0T(\okg) \,\, \stackrel{i_*}{\longrightarrow} \,\, K_0(\okg).\]

{\bf Proposition 3.1.} Let $\gamma$ be one of the 
operations $\sigma^k$, $\lambda^k$,
$\psi^k$ defined in section 1. Then the following diagram commutes:
\[\begin{array}{ccc}
K_1(K_\wp \Gammam) & \stackrel{\Phi}{\longrightarrow} & K_0(\okg) \\
\\
\downarrow \gamma && \downarrow \gamma\\
\\
K_1(K_\wp \Gammam) & \stackrel{\Phi}{\longrightarrow} & K_0(\okg). 
\end{array}\]

We will give an algebraic proof and a topological proof of this proposition. The 
algebraic proof is elementary and explicit, but consists of lengthy
computations. The topological proof is conceptual and uses higher $K$-theory.
We start with a lemma which will be used in the algebraic proof.

{\bf Lemma 3.2.} Let $F$ be a f.~g.\ free $\okg$-module and $\alpha \in 
{\rm GL}(F\otimes K_\wp)$. We write $\alpha = \beta/ \pi$ with
$\beta \in \End(F\otimes \cO_{K_\wp}) \cap {\rm GL}(F\otimes K_\wp)$
and $\pi \in \cO_K\backslash\{0\}$
and we put 
\[C(\beta) := \ker(F \hookrightarrow F\otimes \cO_{K_\wp} \rightepi \,\,
\coker(\beta)).\]
Then we have:\\
(a) $\Phi((F\otimes K_\wp, \alpha)) = [C(\beta)] - [F]$ in 
$K_0(\okg)$. \\
(b) For all $k_1, \ldots, k_r \in \NN$ such that $\gcd(k_1,\ldots, k_r, 
\ord(\Gammam)) =1$ and for all $\beta_1, \ldots, \beta_r \in \End(F\otimes
\cO_{K_\wp}) \cap {\rm GL}(F \otimes K_\wp)$, the canonical
$\okg$-linear homomorphism
\[ \Sym^{k_1}(C(\beta_1)) \otimes \cdots \otimes \Sym^{k_r}(C(\beta_r))
\ra C((\Sym^{k_1}(\beta_1) \otimes \cdots \otimes \Sym^{k_r}(\beta_r))\]
is bijective; here, $\Sym^{k_1}(\beta_1) \otimes \cdots \otimes 
\Sym^{k_r}(\beta_r)$ is considered as an endomorphism of 
$\Sym^{k_1}(F\otimes \cO_{K_\wp}) \otimes \cdots \otimes 
\Sym^{k_r}(F\otimes \cO_{K_\wp}) \cong \Sym^{k_1}(F) \otimes \cdots 
\otimes \Sym^{k_r}(F) \otimes \cO_{K_\wp}$ which is free over $\okpg$ by Proposition
1.1.

{\bf Proof.}  (a) In the commutative diagram 
\[ \begin{array}{ccccccccl} 
0 & \ra &
C(\beta) & \ra & F & \ra & \coker(\beta) & \ra & 0\\ 
\\ 
&& \downarrow && \downarrow && \| \\ 
\\
0 & \ra & F \otimes \cO_{K_\wp} & \stackrel{\beta}{\ra} &
F\otimes \cO_{K_\wp} & \ra & \coker(\beta) & \ra & 0, \end{array}\] 
the
composition $F \hookrightarrow F\otimes \cO_{K_\wp} \rightepi \,\, 
\coker(\beta)$ is
surjective since $(F\otimes \cO_{K_\wp})/F$ is a divisible Abelian group and
$\coker(\beta)$ is a f.~g.\ torsion group.  Thus, the upper sequence is exact.  Since
$\coker(\beta)$ is cohomologically trivial, it has $\okg$-projective dimension
$1$ (cf.\ Proposition 4.1 in \cite{C}).  Hence, the module $C(\beta)$ is
projective over $\okg$ and we have \[\Phi((F\otimes K_\wp, \beta)) = [C(\beta)] -
[F] \quad {\rm in} \quad K_0(\okg).\] For $\beta = \pi$, we obviously have
$C(\beta) \cong F$; consequently, 
$\Phi((F\otimes K_\wp, \pi)) =0$.  This implies assertion (a).\\ 
(b) We write $\Sym^k(-)$ for $\Sym^{k_1}(-) \otimes \cdots
\otimes \Sym^{k_r}(-)$.  The canonical inclusion $C(\beta) \hookrightarrow
F\otimes \cO_{K_\wp}$ obviously induces an isomorphism $C(\beta) \otimes
\cO_{K_\wp} \, \, \tilde{\ra} \,\, F \otimes \cO_{K_\wp}$.  Thus, the canonical
homomorphism $\Sym^k(C(\beta)) \otimes \cO_{K_\wp} \ra \Sym^k(F \otimes
\cO_{K_\wp})$ is bijective.  Hence, in the commutative diagram with exact rows
\[\begin{array}{ccccccccl} 
0 & \ra & \Sym^k(C(\beta)) & \stackrel{i}{\ra} &
\Sym^k(F) & \ra & \coker(i) & \ra & 0 \\ 
\\ 
&& \downarrow && \| && \downarrow \\
\\ 
0 & \ra & C(\Sym^k(\beta)) & \ra & \Sym^k(F) & \ra & \coker(\Sym^k(\beta)) &
\ra & 0\\ 
\\ 
&& \downarrow && \downarrow && \|\\ 
\\ 0& \ra & \Sym^k(F\otimes
\cO_{K_\wp}) & \stackrel{\Sym^k(\beta)}{\ra} & \Sym^k(F \otimes \cO_{K_\wp}) & \ra
& \coker(\Sym^k(\beta)) & \ra & 0, \end{array}\] 
the right vertical composition
is bijective.  Then, the upper right vertical homomorphism is bijective, too.
Hence, the upper left vertical homomorphism is bijective, as was to be shown.

{\bf Algebraic Proof of Proposition 3.1.}  We use additive notation for
$K_1(K_\wp \Gammam)$.  Let $(F \otimes K_\wp, \beta/\pi)$ be an element of
$K_1(K_\wp \Gammam)$ as in Lemma 3.2.  Then the element 
\[[(F\otimes K_\wp,
\beta/\pi)] - [(F\otimes K_\wp, 1)] \quad {\rm of} \quad 
\tilde{K}_0(\ZZ, K_\wp \Gammam)\] 
is a
preimage of $(F \otimes K_\wp, \beta/\pi)$ under the canonical homomorphism
$\tilde{K}_0(\ZZ, K_\wp \Gammam) \ra K_1(K_\wp \Gammam)$ (see also Theorem 1.6(f)).
Using Theorem 1.6(c)(i) and (f), Lemma 3.2(a), Lemma 3.2(b), 
Theorem 1.6(c)(i),
and Lemma 3.2(a) successively, we obtain:  
\begin{eqnarray*}
\lefteqn{\Phi\left(\sigma^k\left((F\otimes K_\wp, \beta/\pi)\right)\right) = }\\ 
&& = \Phi\Bigg(\sum_{{a, b_1,
\ldots, b_u \ge 1}\atop{a+b_1 + \cdots + b_u=k}} 
(-1)^u\Big(\Sym^a(F\otimes K_\wp)
\otimes \Sym^{b_1}(F\otimes K_\wp) \otimes \cdots \otimes \Sym^{b_u}(F \otimes
K_\wp), \\
&& \hspace*{24em}\Sym^a(\beta)/\pi^a \otimes 1 \otimes \cdots \otimes 1\Big)\Bigg)\\ 
&& = \sum_{{a,
b_1, \ldots, b_u \ge 1}\atop{a+b_1 + \cdots + b_u=k}} (-1)^u\Big([C(\Sym^a(\beta)
\otimes 1 \otimes \cdots \otimes 1)] - \\
&& \hspace*{12em} [\Sym^a(F) \otimes \Sym^{b_1}(F) \otimes
\cdots \otimes \Sym^{b_u}(F)]\Big)\\ 
&& = \sum_{{a, b_1, \ldots, b_u \ge 1}\atop{a+b_1
+ \cdots + b_u=k}} (-1)^u\Big([\Sym^a(C(\beta)) \otimes \Sym^{b_1}(F) \otimes \cdots
\otimes \Sym^{b_u}(F)] - \\
&& \hspace*{15em} [\Sym^a(F) \otimes \Sym^{b_1}(F) \otimes \cdots \otimes
\Sym^{b_u}(F)]\Big)\\ 
&& = \sum_{{a \ge 0, b_1, \ldots, b_u \ge 1}\atop{a+b_1 + \cdots
+ b_u=k}} (-1)^u [\Sym^a(C(\beta)) \otimes \Sym^{b_1}(F) \otimes \cdots \otimes
\Sym^{b_u}(F)] \\ 
&& = \sigma^k \left([C(\beta)] - [F]\right)\\ 
&& = \sigma^k
\left(\Phi\left((F\otimes K_\wp, \beta/\pi)\right)\right).  
\end{eqnarray*} 
This proves Proposition
3.1 for $\gamma = \sigma^k$.  For $\gamma = \lambda^k$ or $\gamma = \psi^k$,
Proposition 3.1 can be shown analogously or can be deduced from the assertion for
$\gamma = \sigma^k$ using Theorem 1.6(c) and (d).

{\bf Topological Proof of Proposition 3.1.} Let $K_0T_{\rm pr}(\okpg)$
denote the Grothendieck group of all torsion modules over $\okpg$ which
have a $\okpg$-projective resolution of length $1$. Then the canonical
homomorphism $K_0T(\okpg) \ra K_0T_{\rm pr}(\okpg)$ is bijective. The
injectivity has been proved in Theorem 1(ii) on p.\ 3 in \cite{F2}. The
surjectivity follows from the fact that, in the commutative diagram
\[\begin{array}{ccccccc}
K_1(K_\wp \Gammam) & \stackrel{\partial}{\ra} & K_0T(\okpg)\\
\\
\| &&\downarrow\\
\\
K_1(K_\wp \Gammam) & \stackrel{\partial}{\ra} & K_0T_{\rm pr}(\okpg) &
\stackrel{i_*}{\ra} & K_0(\okpg) & \stackrel{\rm can}{\ra} &
K_0(K_\wp\Gammam),
\end{array}\]
the lower sequence is exact (cf.\ Theorem 1(ii) on p.\ 3 
in \cite{F2}) and that the
canonical homomorphism $K_0(\okpg) \ra K_0(K_\wp \Gammam)$ is injective
(cf.\ Theorem 6.1 in \cite{Sw}). In the section ``The localization
theorem for projective modules'' in \cite{Gr1}, it is shown that
$K_0T_{\rm pr}(\okpg)$ is isomorphic to the group of connected components
of the homotopy fibre of the canonical continuous map from the
$K$-theory space associated with $\okpg$ to the $K$-theory space
associated with $K_\wp \Gammam$. Thus, the continuous map $\sigma^k$
on the $K$-theory space associated with $\okpg$ defined in section 1
induces a symmetric power operation $\sigma^k$ on $K_0T_{\rm pr}(\okpg)$
which in particular commutes with the connecting homomorphism 
$\partial: K_1(K_\wp \Gammam) \ra K_0T_{\rm pr}(\okpg)$. Analogously,
the Grothendieck group $K_0T_\wp(\okg)$ of all $\wp$-torsion modules
over $\okg$ of $\okg$-projective dimension $1$ carries a symmetric
power operation $\sigma^k$ which commutes with the localization
isomorphism $K_0T_{\wp}(\okg) \, \, \tilde{\ra} \,\, K_0T_{\rm pr}(\okpg)$
and with the resolution map $i_*:K_0T_\wp(\okg) \ra K_0(\okg)$. Since $\Phi$
coincides with the composition
\[K_1(K_\wp \Gammam) \,\, \stackrel{\partial}{\ra} \,\,
K_0T_{\rm pr}(\okpg) 
\,\, \tilde{\leftarrow} \,\, K_0T_\wp(\okg) \,\, 
\stackrel{i_*}{\ra} \,\, K_0(\okg)\]
(cf.\ Proposition 2.5 on p.\ 183 in \cite{Sh}), 
Proposition 3.1 is proved for $\sigma^k$. The assertions 
for $\lambda^k$ and $\psi^k$ can be deduced from this as in the algebraic proof. 

{\bf Remark 3.3.} Let $X$ be a scheme and $U$ an open subscheme of $X$.
As in the proof above, one can define power operations on the 
homotopy fibre of the canonical continuous map from the $K$-theory
space of $X$ to the $K$-theory space of $U$. If $U$ is affine and
$X\backslash U$ is a divisor, it is shown in the section ``The
localization theorem \ldots'' in \cite{Gr1} that this
homotopy fibre is homotopy equivalent to the $K$-theory space associated
with the category $\cal  H$ consisting of all quasi-coherent sheaves on 
$X$ which are zero on $U$ and which have a resolution by vector bundles
on $X$ of length $1$. It is likely that one can use the ideas developed
in the algebraic proof to describe these operations on the $K$-theory
space of $\cal H$ simplicially, and, finally, on $K_0({\cal H})$ algebraically,
i.e., in terms of modules. I hope to say more on this in a future paper.

Let $\cM_1, \cM_2, \ldots$ be the categories defined in section 1. 

{\bf Proposition 3.4.} For any $1 \le i \le k-1$, the composition
\[{\rm Cl}(\okg) \times {\rm Cl}(\okg) 
%\subset K_0(\okg) \times K_0(\okg)
\,\, \stackrel{\sigma^i \times \sigma^{k-i}}{\longrightarrow} \,\,
K_0(\cM_i) \times K_0(\cM_{k-i}) \,\, \stackrel{\otimes}{\rightarrow} \,\,
K_0(\okg) \]
is the zero map. In particular, the restriction of $\sigma^k: 
K_0(\okg) \ra K_0(\okg)$ to the locally free classgroup ${\rm Cl}(\okg)$
is a homomorphism. Analogous assertions hold for exterior power operations,
too.

{\bf Proof.} Let $(x,y)\in {\rm Cl}(\okg) \times {\rm Cl}(\okg)$. Let 
$P_1, P_2, Q_1, Q_2$ be f.~g.\ projective $\okg$-modules such that
$x= [P_1] - [Q_1]$ and $y = [P_2] - [Q_2]$. By Theorem A in \cite{Sw},
we may assume that $P_1, P_2, Q_1, Q_2$ are projective left ideals in
$\okg$. Then it suffices by Theorem 1.6(c)(i) to show Proposition 3.4 for 
elements $x,y$ of the form $x=[\okg] - [I]$ and $y= [\okg] - [J]$ where
$I,J$ are projective left ideals in $\okg$. Furthermore, we may assume
by Theorem A in \cite{Sw} that the annihilators of the torsion modules
$\okg/I$ and $\okg/J$ are relatively prime. By Theorem 1.6(c)(i), we have
in $K_0(\cM_i)$:
\begin{eqnarray*}
\lefteqn{\sigma^i(x) = \sum_{{a \ge 0, b_1, \ldots, b_u \ge 1}\atop
{a+ b_1 + \cdots + b_u = i}} (-1)^u 
\left[\Sym^a(\okg) \otimes \Sym^{b_1}(I) \otimes \cdots \otimes \Sym^{b_u}(I)
\right] = }\\
&& = \sum_{{a , b_1, \ldots, b_u \ge 1}\atop{a+ b_1 + \cdots + b_u = i}} 
(-1)^u \bigg(\left[\Sym^a(\okg) \otimes \Sym^{b_1}(I) \otimes \cdots \otimes 
\Sym^{b_u}(I)\right] - \\
&& \hspace*{15em} \left[\Sym^a(I)\otimes \Sym^{b_1}(I) \otimes \cdots \otimes 
\Sym^{b_u}(I)\right]\bigg).
\end{eqnarray*}
For any $a\ge 1$ and for any prime ideal $\wp \not= 0$ in $\cO_K$ which 
does not lie in the support of $\okg/I$, the localization of the
canonical homomorphism 
\[\Sym^a(I) \ra \Sym^a(\okg)\]
at $\wp$ is bijective. In particular, $\Sym^a(I) \ra \Sym^a(\okg)$ is 
injective and the support of the torsion module $\Sym^a(\okg)/\Sym^a(I)$
is contained in the support of the torsion module $\okg/I$. The same holds
for $y$ and $J$. Consequently,
it suffices to show the following assertion: Let $M_1 \subseteq N_1$ and
$M_2 \subseteq N_2$ be modules in $\cM_i$ and $\cM_{k-i}$, respectively,
such that $N_1/M_1$ and $N_2/M_2$ are torsion modules whose 
annihilators are relatively prime. Then we have:
\[([N_1] - [M_1]) \cdot ([N_2] - [M_2]) =0 \quad {\rm in} \quad
K_0(\okg).\]
This immediately follows from the fact that, in the commutative diagram
with exact rows
\[\begin{array}{ccccccccl}
0 & \ra & M_1 \otimes M_2 & \ra & N_1 \otimes M_2 & \ra & N_1/M_1 \otimes M_2
& \ra & 0 \\
\\
&& \downarrow && \downarrow && \downarrow\\
\\
0 & \ra & M_1 \otimes N_2 & \ra & N_1 \otimes N_2 & \ra & N_1/M_1 \otimes N_2
& \ra & 0.
\end{array}\]
the right vertical homomorphism is bijective. \\
The equality $\sigma^k(x+y)
= \sum_{i=0}^k \sigma^i(x) \cdot \sigma^{k-i}(y)$ (for $x, y \in {\rm Cl}(\okg)$)
now shows that $\sigma^k: {\rm Cl}(\okg) \ra {\rm Cl}(\okg)$ is a
homomorphism.

As for Proposition 3.1, we give a topological proof for Proposition 3.4
with a more conceptual argument.

{\bf Topological proof of Proposition 3.4 (Sketch).} We will consider
only the case $k=2$, $i=1$; i.e., we will only show
that the multiplication
\[{\rm Cl}(\okg) \times {\rm Cl}(\okg) \ra K_0(\okg)\]
is trivial. The general statement can then be proved by generalizing the following
argument as in the algebraic proof.\\
As in the topological proof of Proposition 3.1, the tensor product
induces a multiplication on the various $K$-theory spaces and, in particular,
on the $K$-groups $K_1(K_\wp \Gammam)$, $K_0T(\okg)$, $K_0(\okg)$, \ldots\ 
which are compatible with the maps in the local or global localization
sequences. On $K_1(K_\wp \Gammam)$ (as on all higher $K$-groups), this
multiplication is linear. Hence, it is trivial. Since the connecting 
homomorphism $\partial: K_1(K_\wp \Gammam) \ra K_0T(\okpg)$ is surjective,
the multiplication on $K_0T(\okpg)$ is trivial, too. As in the algebraic 
proof, one shows furthermore that for all prime ideals $\wp \not= \wq$ in 
$\cO_K$, the multiplication $K_0T_\wp(\okg) \times K_0T_\wq(\okg) \ra
K_0(\okg)$ is trivial. Therefore, the multiplication $K_0T(\okg) \times
K_0T(\okg) \ra K_0(\okg)$ is trivial. This implies that the 
multiplication on ${\rm Cl}(\okg)$ is trivial, as was to be shown.

Now let $\gamma$ be one of the operations $\sigma^k$, $\lambda^k$, or
$\psi^k$ defined in section 1. For any prime ideal 
$\wp \not= 0$ in $\cO_K$, the operation
$\gamma$ commutes with the canonical homomorphism $K_1(\okpg) \ra
K_1(K_\wp \Gammam)$ by Theorem 1.6(b) and hence induces an operation
$\gamma$ on $K_0T(\okpg)$. Let $\gamma$ also denote the direct sum
of these operations on $K_0T(\okg) \cong \oplus_\wp K_0T(\okpg)$.
Then, $\gamma$ commutes with the connecting homomorphism
$\partial: K_1(K\Gammam) \ra K_0T(\okg)$ by Theorem 1.6(b). Consequently,
$\gamma$ induces an operation on ${\rm Cl}(\okg) \cong
\coker(\partial)$ denoted by $\gamma$ again. On the other hand, the operation $\gamma$ on 
$K_0(\okg)$ obviously induces an operation $\gamma$
on ${\rm Cl}(\okg)$ by restricting. 

{\bf Corollary 3.5.} These two definitions of the operation $\gamma$ on
${\rm Cl}(\okg)$ coincide. 

{\bf Proof.} This immediately follows from Propositions 3.1 and 3.4.

For any prime ideal $\wp \not= 0$ in $\cO_K$, let $G_\wp := 
\Gal(\Kpq/K_\wp)$ denote the absolute Galois group. By 
Chapter 1 in \cite{T}, the pairing
\[\begin{array}{ccl}
K_0(\Kpqg) \times K_1(\Kpqg) & \ra & \Kpq^\times \\
([P], (M,\alpha)) & \mapsto & \det_\Kpq 
(\Hom_\Kpqg (P,\alpha)|_{\Hom_\Kpqg(P,M)})
\end{array}\]
(here, $\Hom_\Kpqg (P,\alpha)$ denotes the automorphism of $\Hom_\Kpqg(P,M)$
induced by $\alpha$) induces an isomorphism
\[K_1(K_\wp \Gammam) \,\, \tilde{\ra} \,\, 
\Hom_{G_\wp} (K_0(\Kpqg), \Kpq^\times).\]
Note that the determinant map used here coincides with the determinant map
used in \cite{T} (see Proposition 2.3 in \cite{Que}). We identify the 
Grothendieck group $K_0(\Kpqg)$ with the classical ring of virtual characters
of $\Gammam$. Let $\langle-,-\rangle: K_0(\Kpqg) \times K_0(\Kpqg) \ra
\ZZ$ denote the classical character pairing.

Now let $k$ be an arbitrary natural number. Let $\hat{\psi}^k
\in \End_{G_\wp}(K_0(\Kpqg))$ denote 
the adjoint of the Adams operation $\psi^k$ with respect to $\langle-,-\rangle$.

{\bf Lemma 3.6.} (a) For all characters $\chi \in K_0(\Kpqg)$ and
all $\gamma \in \Gammam$, we have:
\[\hat{\psi}^k(\chi)(\gamma) = \sum_{\tau^k=\gamma} \chi (\tau) \quad
{\rm in} \quad \Kpq.\]
In particular, we have $\hat{\psi}^k = \psi^{k'}$ if 
$\gcd(k, \ord(\Gammam)) =1$ and $kk' \equiv 1$ mod $e(\Gammam)$. \\
(b) Let $\Gammam$ be an Abelian group. Let $\phi: \Kpqg \ra \Kpqg$ denote
the ring homomorphism given by $\gamma \mapsto \gamma^k$. Then we have for
all f.~g.\ $\Kpqg$-modules $M$:
\[\hat{\psi}^k([M]) = \left[(\Kpqg)_\phi \otimes_\Kpqg M\right] \quad
{\rm in} \quad K_0(\Kpqg).\]

{\bf Proof.} (a) Let the class function $\tilde{\psi}^k(\chi): \Gammam \ra 
\CC$ be given by $\tilde{\psi}^k(\chi)(\gamma) = 
\sum_{\tau^k=\gamma} \chi(\tau)$
for $\gamma \in \Gammam$. Then we have for all characters $\theta \in 
K_0(\CC\Gammam)$:
\begin{eqnarray*}
\lefteqn{\langle \tilde{\psi}^k(\chi), \theta\rangle = \frac{1}{\ord(\Gammam)}
\sum_{\gamma \in \Gammam} 
\tilde{\psi}^k(\chi)(\gamma)\overline{\theta(\gamma)}}\\
&& = \frac{1}{\ord(\Gammam)} \sum_{\gamma \in \Gammam} \sum_{\tau^k = \gamma}
\chi(\tau) \overline{\theta(\gamma)}
= \frac{1}{\ord(\Gammam)} \sum_{\tau \in \Gammam} \chi(\tau) 
\overline{\theta(\tau^k)}\\
&& = \frac{1}{\ord(\Gammam)} \sum_{\tau \in \Gammam} \chi(\tau)
\overline{\psi^k(\theta)(\tau)}
= \langle \chi, \psi^k(\theta)\rangle = \langle \hat{\psi}^k(\chi), \theta
\rangle.
\end{eqnarray*}
(Here, we have used Proposition (12.8) on p.\ 317 in \cite{CR}.) This shows
the main assertion of Lemma 3.6(a). 
If $\gcd(k,\ord(\Gammam)) =1$, then we have $\{\tau \in \Gammam:
\tau^k = \gamma\} = \{\gamma^{k'}\}$ for all $\gamma \in \Gammam$. This
proves $\hat{\psi}^k = \psi^{k'}$ in this case.\\
(b) For all $\Kpqg$-modules $N$, we have:
\begin{eqnarray*}
\lefteqn{ \langle (\Kpqg)_\phi \otimes_\Kpqg M, N \rangle =
\dim_\Kpq \Hom_\Kpqg((\Kpqg)_\phi \otimes_\Kpqg M, N) =} \\
&& = \dim_\Kpq \Hom_\Kpqg(M,\psi^k(N))
 = \langle M, \psi^k(N)\rangle = \langle \hat{\psi}^k(M), N\rangle.
\end{eqnarray*}
This proves Lemma 3.6(b).

We now assume that the prime divisors of $\ord(\Gammam)$ are unramified in
$K$. Furthermore, let $k$ be odd if $\Gammam$ has irreducible symplectic 
characters. Under these assumptions, Ph.\ Cassou-Nogu\`es and M.\ Taylor
have shown that the endomorphisms $\psi_k^{\rm CNT} := \Hom(\psi^k, \Kpq^\times)$
of $K_1(\Kpqg) \cong \Hom_{G_\wp}(K_0(\Kpqg), \Kpq^\times)$, $\wp \in 
\Max(\cO_K)$, induce an endomorphism of ${\rm Cl}(\okg)$
denoted by $\psi_k^{\rm CNT}$ again (cf.\ Chapter 9 in \cite{T}). In the
paper \cite{BC}, D.\ Burns and T.\ Chinburg have raised the question whether 
there is an algebraic interpretation of this endomorphism. The 
following theorem gives an affirmative answer to 
this question in the case $\gcd(k,\ord(\Gammam))=1$. For this, let 
$\sigma^k$ and $\psi^k$ denote (the restriction of) the symmetric power
operation and $k$-th Adams operation on ${\rm Cl}(\okg) \subset K_0(\okg)$
constructed in section 1. 

{\bf Theorem 3.7.} Let $\gcd(k, \ord(\Gammam)) =1$ and $k'\in \NN$ such
that $kk' \equiv 1$ mod $e(\Gammam)$. Then we have:
\[\psi_{k'}^{\rm CNT} = \sigma^k \quad {\rm on} \quad {\rm Cl}(\okg)\]
and
\[ k\cdot \psi_{k'}^{\rm CNT} \oplus \id = \psi^k \quad {\rm on} \quad
{\rm Cl}(\okg) \oplus \ZZ[\okg] = K_0(\okg).\]

{\bf Proof.} By Theorem 3.3 in \cite{KoAdHi}, the operation $\lambda^k$
on $K_1(\Kpqg)$ corresponds to the endomorphism $(-1)^{k-1} \Hom(\hat{\psi}^k,
\Kpq^\times)$ of $\Hom_{G_\wp}(K_0(\Kpqg), \Kpq^\times)$ via the
isomorphism 
\[K_1(\Kpqg) \cong \Hom_{G_\wp}(K_0(\Kpqg), \Kpq^\times).\]
Thus, the operation $\sigma^k$ corresponds to the endomorphism
$\Hom(\psi^{k'}, \Kpq^\times)$ by Theorem 1.6(d)(ii) and Lemma 3.6(a).
By construction (see Chapter 9 in \cite{T}), $\psi_{k'}^{\rm CNT}$ on
${\rm Cl}(\okg)$ is induced by $\Hom(\psi^{k'}, \Kpq^\times)$ on 
$K_1(\Kpqg)$, $\wp \in \Max(\cO_K)$. Now the equality $\psi_{k'}^{\rm CNT}
= \sigma^k$ on ${\rm Cl}(\okg)$ follows from Corollary 3.5. The equality
$k \cdot \psi_{k'}^{\rm CNT} = \psi^k$ on ${\rm Cl}(\okg)$ can be shown
analogously. Finally, the equality $\psi^k = \id$ on $\ZZ[\okg] \subseteq
K_0(\okg)$ has already been shown in Theorem 1.6(e).

{\bf Corollary 3.8.} We have $\sigma^k = \sigma^{k+e(\Gammam)}$ on
${\rm Cl}(\okg)$. 

{\bf Proof.} This follows from Theorem 1.3(2) on p.\ 99 in \cite{T} and
Theorem 3.7.

The following corollary answers the question (1.12) in Chapter 9 of \cite{T}.

{\bf Corollary 3.9.} Let $c: {\rm Cl}(\okg) \hookrightarrow 
K_0(\okg) \ra K_0(\Gammam, \cO_K)$ denote the Cartan homomorphism. Then
the following diagrams commute:
\[\begin{array}{ccc}
{\rm Cl}(\okg) & \stackrel{c}{\ra} & K_0(\Gammam, \cO_K)\\
\\
\phantom{\psi_{k'}^{\rm CNT}} \downarrow \psi_{k'}^{\rm CNT} && 
\phantom{\sigma^k} \downarrow \sigma^k\\
\\
{\rm Cl}(\okg) & \stackrel{c}{\ra} & K_0(\Gammam, \cO_K)
\end{array}\qquad
\begin{array}{ccc}
{\rm Cl}(\okg) \oplus \ZZ[\okg]& \stackrel{c}{\ra} & K_0(\Gammam, \cO_K)\\
\\
\phantom{k \cdot \psi_{k'}^{\rm CNT} \oplus \id } \downarrow 
k \cdot \psi_{k'}^{\rm CNT} \oplus \id && 
\phantom{\psi^k} \downarrow \psi^k\\
\\
{\rm Cl}(\okg) \oplus \ZZ[\okg] & \stackrel{c}{\ra} & K_0(\Gammam, \cO_K).
\end{array}\]

{\bf Proof.} This immediately follows from Theorem 1.6(a) and Theorem 3.7.
In addition to the algebraic or topological arguments used in the proof
of Theorem 3.7, we give a geometric proof of this corollary since 
Proposition 3.10 which will be used in this proof is of independent 
interest.\\
By Theorem 1.3(3) on p.\ 99 in \cite{T}, the following diagram commutes 
for any prime ideal $\wp \not= 0 $ in $\cO_K$:
\[\begin{array}{ccccc}
K_0T(\okpg) & \stackrel{c}{\longrightarrow} & K_0T(\Gammam, \cO_{K_\wp}) & 
\cong & K_0(\Gammam, \cO_K/\wp) \\
\\
\phantom{\psi_{k'}^{\rm CNT}} \downarrow \psi_{k'}^{\rm CNT} &&&&
\phantom{\psi^k} \downarrow \psi^k\\
\\
K_0T(\okpg) & \stackrel{c}{\longrightarrow} & K_0T(\Gammam, \cO_{K_\wp}) & 
\cong & K_0(\Gammam, \cO_K/\wp). 
\end{array}\]
Now Corollary 3.9 follows from the following Proposition 3.10 (and
Theorem 1.6(e)) where, for the first diagram in Corollary 3.9, we
in addition use the fact that, for all prime ideals $\wp \not= \wq$ in
$\cO_K$ and all $x\in K_0(\Gammam, \cO_K/\wp)$, $y \in K_0(\Gammam, \cO_K/\wq)$,
the product $i_*(x) \cdot i_*(y)$ vanishes in $K_0(\Gammam, \cO_K)$ (cf.\
the end of the (algebraic) proof of Proposition 3.4).

{\bf Proposition 3.10.} For any prime ideal $\wp \not= 0$ in $\cO_K$ and
any $j\ge 0$, the following diagrams commute:
\[\begin{array}{ccc}
K_0(\Gammam, \cO_K/\wp) & \stackrel{i_*}{\rightarrow} & K_0(\Gammam, \cO_K)\\
\\
\phantom{\psi^j} \downarrow \psi^j && \phantom{\sigma^j} \downarrow \sigma^j\\
\\
K_0(\Gammam, \cO_K/\wp) & \stackrel{i_*}{\rightarrow} & K_0(\Gammam, \cO_K)
\end{array}\qquad \qquad
\begin{array}{ccc}
K_0(\Gammam, \cO_K/\wp) & \stackrel{i_*}{\rightarrow} & K_0(\Gammam, \cO_K)\\
\\
\phantom{j \cdot \psi^j} \downarrow j \cdot \psi^j && 
\phantom{\psi^j} \downarrow \psi^j\\
\\
K_0(\Gammam, \cO_K/\wp) & \stackrel{i_*}{\rightarrow} & K_0(\Gammam, \cO_K).
\end{array}\]
Here, $i_*$ denotes the composition of the canonical map 
$K_0(\Gammam, \cO_K/\wp) \ra K'_0(\okg)$ with the isomorphism
$K'_0(\okg) \cong K_0(\Gammam, \cO_K)$ (cf.\ Notations). 

{\bf Proof.} The conormal module $\wp/\wp^2$ of the closed immersion
$i: \Spec(\cO_K/\wp) \hookrightarrow \Spec(\cO_K)$ is trivial. Thus,
the Bott multiplier $\theta^j(i) := 1 + [\wp/\wp^2] + \cdots + 
[(\wp/\wp^2)^{j-1}]$ associated with $i$ equals $j$ in $K_0(\Gammam, \cO_K/\wp)$. Now,
the commutativity of the second diagram follows from the equivariant
Adams-Riemann-Roch theorem without denominators (see Corollary (5.2) in
\cite{KoARR}). The commutativity of the first diagram follows from the
(more general) 
equivariant Riemann-Roch theorem (see Satz (5.1) in \cite{KoARR}) 
applied to the natural operation $\sigma^j$. For this, we have to show
\[\sigma^j(\wp/\wp^2, x) = \psi^j(x) \quad {\rm in} \quad 
K_0(\Gammam, \cO_K/\wp)\]
for all $x\in K_0(\Gammam, \cO_K/\wp)$. Here, $\sigma^j(\wp/\wp^2, x)$ is
defined as follows (see section 5 in \cite{KoARR}): Let the
polynomial ring $K_0(\Gammam, \cO_K/\wp)[N]$ be equipped with the unique (special)
$\lambda$-ring structure such that $K_0(\Gammam, \cO_K/\wp)$ is a 
$\lambda$-subring and such that $N$ is of $\lambda$-degree $1$; then, 
$\sigma^j(\wp/\wp^2, x)$ is defined to be the value of the polynomial
$\sigma^j(x \cdot (1-N))/(1-N)$ at the place $N=1$. As for exterior powers
(see p.\ 5 in \cite{FL}), there is
a universal polynomial $Q_j \in \ZZ[X_1, \ldots, X_j; Y_1, \ldots, Y_j]$
such that 
\[\sigma^j(x \cdot (1-N)) = Q_j(\sigma^1(x), \ldots, \sigma^j(x);
\sigma^1(1-N), \ldots, \sigma^j(1-N)).\]
Since $\sigma^i(1-N) = 1-N$ for
all $i\ge 0$ and since $Q_j$ is homogeneous of weight $j$ in both sets 
of variables, we have:
\[\sigma^j(\wp/\wp^2, x) = Q_j(\sigma^1(x), \ldots, \sigma^j(x); 0, \ldots,
0,1) \quad {\rm in} \quad K_0(\Gammam, \cO_K/\wp).\]
If $L_1, \ldots, L_n$ are of modules of rank 1, $M:= L_1 \oplus \cdots 
\oplus L_n$, and $N$ is an arbitrary module, we have:
\begin{eqnarray*}
\lefteqn{Q_j\left([\Sym^1(M)], \ldots, [\Sym^j(M)]; 
[\Sym^1(N)], \ldots, [\Sym^j(N)]\right) =}\\
&& = [\Sym^j((L_1 \oplus \cdots \oplus L_n) \otimes N)] \\
&& = [\Sym^j(L_1 \otimes N \oplus \cdots \oplus L_n\otimes N)]\\
&& = \sum_{b_1+ \cdots + b_n = j} [\Sym^{b_1}(L_1\otimes N)\otimes \cdots
\otimes \Sym^{b_n}(L_n\otimes N)]\\
&& = \sum_{b_1+ \cdots + b_n = j} [L_1^{\otimes b_1} \otimes \cdots \otimes
L_n^{\otimes b_n} \otimes \Sym^{b_1}(N) \otimes \cdots \otimes \Sym^{b_n}(N)].
\end{eqnarray*}
Thus, for $x= \sum_{i=1}^n l_i$ with elements $l_i$ of $\lambda$-degree $1$, we
have:
\[Q_j(\sigma^1(x), \ldots, \sigma^j(x); Y_1, \ldots, Y_j) = 
\sum_{b_1+\cdots + b_n=j} l_1^{b_1}   \cdots   l_n^{b_n}\cdot
Y_{b_1}   \cdots   Y_{b_n}\]
(with $Y_i := 0$ for $i > j$) and, hence,
\[Q_j(\sigma^1(x), \ldots, \sigma^j(x); 0, \ldots, 0, 1) = l_1^j + \cdots
+ l_n^j = \psi^j(x).\]
Using the splitting principle (see section (2.5) in \cite{KoARR}), we
finally obtain the equality $\sigma^j(\wp/\wp^2, x) = \psi^j(x)$ for 
an arbitrary $x \in K_0(\Gammam, \cO_K/\wp)$. This completes the proof of 
Proposition 3.10.

{\bf Remark 3.11.} The construction
of $\psi_k^{\rm CNT}$ in \cite{T} only shows that, for any
prime ideal $\wp \not= 0 $ in $\cO_K$, the image of the canonical 
homomorphism 
\[K_1(\okpg) \ra K_1(K_\wp \Gammam) \cong \Hom_{G_\wp}
(K_0(\Kpqg), \Kpq^\times)\]
is invariant under the endomorphism $\Hom(\psi^k, \Kpq^\times)$. If  
$p:= \char(\cO_K/\wp)$ does not divide $\gcd(k,\ord(\Gammam))$, 
the considerations above moreover show that there
exists an endomorphism of $K_1(\okpg)$ which is compatible with
$\Hom(\psi^k, \Kpq^\times)$ on $K_1(K_\wp \Gammam)$ and that the 
assumption ``$K_\wp$ is absolutely unramified'' used in \cite{T} is not
necessary for this. There remains open the question whether, in the
case $p \, | \, \gcd(k,\ord(\Gammam))$, there exists an endomorphism
of $K_1(\okpg)$ which is compatible with $\Hom(\psi^k, \Kpq^\times)$. This
question has an affirmative answer in the following two (extreme) cases:\\
(a) Let $\Gammam$ be an Abelian group. Let $\phi: \okpg \ra \okpg$ denote 
the ring homomorphism given by $\gamma \mapsto \gamma^k$. Then $\phi$
induces an endomorphism on $K_1(\okpg)$ which is compatible with 
$\Hom(\psi^k,\Kpq^\times)$ on $K_1(K_\wp \Gammam)$. This is easy to prove or 
follows from
formula (3.18) on p.\ 113 in \cite{T} which holds for arbitrary Abelian
groups and not only for $p$-groups.\\
(b) Let $k$ be equal to the exponent of $\Gammam$. Let 
$\epsilon: \okpg \ra \okpg$, $\sum a_\gamma[\gamma] \mapsto 
\sum a_\gamma [1]$, denote the augmentation homomorphism. Then $\epsilon$
induces an endomorphism on $K_1(\okpg)$ which is compatible with
$\Hom(\psi^k,\Kpq^\times)$ on $K_1(K_\wp \Gammam)$, as one can easily show.

{\bf Remark 3.12.} In the definition of $\psi_k^{\rm CNT}$, it has been
assumed that $k$ is odd if $\Gammam$ has irreducible symplectic characters.
This assumption is obviously implied by the assumption $\gcd(k,\ord(\Gammam))
=1$. Thus, the description of $\psi_k^{\rm CNT}$ given in this section
does not yield any improvement into this direction. The following 
considerations about $K_1(\RR\Gammam)$ will make clear why this assumption
is necessary.\\
Let $j$ denote complex conjugation and let $K_0(\CC\Gammam) = K_0^\RR \oplus
K_0^\CC \oplus K_0^\HH$ be the decomposition of the classical ring
$K_0(\CC \Gammam)$ of virtual characters introduced in I, section 12.6 in 
\cite{Se1}. Then the subgroup $K_0^{\rm s}$ of $K_0(\CC \Gammam)$ 
generated by symplectic characters has the following description:
\[K_0^{\rm s} = 2 K_0^\RR \oplus (1+j) K_0^\CC \oplus K_0^\HH =
(1+j) K_0(\CC\Gammam) + K_0^\HH.\]
Furthermore, we have a canonical isomorphism
\[K_1(\RR\Gammam) \cong \Hom_j^+(K_0(\CC\Gammam), \CC^\times)\]
where $\Hom_j^+(K_0(\CC\Gammam), \CC^\times)$ denotes the group of all
Galois-invariant homomorphisms $f: K_0(\CC\Gammam) \ra \CC^\times$ with
$f(K_0^{\rm s}) \subseteq \RR_+^\times$ (see Chapter 1 in \cite{T}). 
Now, the assertion that 
a Galois-invariant operation $\gamma$ on $K_0(\CC\Gammam)$ induces an
operation on $K_1(\RR\Gammam)$ is equivalent to the assertion that
$\gamma$ maps the subgroup $K_0^\HH$ into $K_0^{\rm s}$. If $\Gammam$
is the quaternion group $Q_8$ and $\gamma = \psi^2$, the latter 
assertion is not true as one can easily show; i.e., the analogue of
$\psi_2^{\rm CNT}$ cannot be defined on $K_1(\RR Q_8)$ (and $K_1(KQ_8)$). 
However, we have 
\[\hat{\psi}^k(K_0^\HH) \subseteq K_0^{\rm s}\]
for all $k\ge 1$. This follows from the fact that $\Hom(\hat{\psi}^k, 
\CC^\times)$ corresponds to the operation $\sigma^k$ on $K_1(\RR\Gammam)$
(see the proof of Theorem 3.7). Alternatively, this can be proved in
the following elementary way: Let $S \in K_0^\HH$ be an irreducible
$\CC\Gammam$-module. Then we have to show that, for all irreducible 
$V\in K_0^\RR$, the number $\langle V, \hat{\psi}^k(S) \rangle$ is even
and that, for all irreducible $V\in K_0^\CC$, the equality
$\langle V, \hat{\psi}^k(S) \rangle = \langle j(V), \hat{\psi}^k(S) \rangle$
holds. The first statement follows from the fact that, for all 
f.~g.\ $\RR\Gammam$-modules $W$, the dimension $\langle \CC \otimes W, S 
\rangle = \dim_\CC\Hom_{\CC\Gammam} (\CC\otimes W, S) $ is even since
$\Hom_{\RR\Gammam}(W,S)$ is a vector space over $\HH$. The second statement
is clear since $j$ fixes $S$.

{\bf Remark 3.13.} Let $\Gammam$ be an abelian group, $\wp \in \Max(\cO_K)$,
and $\phi_k: K_\wp \Gammam \ra K_\wp \Gammam$ given by $\gamma \mapsto 
\gamma^k$. Pulling back the module structure along $\phi_k$ defines an
endomorphism $\phi_k^*$ on $K_1(\Gammam, K_\wp) = K_1(K_\wp \Gammam)$. Then
we have 
\[\phi_k^* = \sigma^k \quad {\rm on} \quad K_1(K_\wp \Gammam)\]
for all $k \ge 1$. 
This follows from Lemma 3.6(b) and the fact that the operation $\sigma^k$
corresponds to $\Hom(\hat{\psi}^k, \Kpq^\times)$ via the isomorphism
$K_1(K_\wp \Gammam) \cong \Hom_{G_\wp}(K_0(\Kpqg), \Kpq^\times)$.
Note that $\phi_k^*$ cannot be defined on $K_1(\okpg)$ in general since
a projective $\okpg$-module need not to be projective again after pulling
back the module structure along $\phi_k$. Thus, it is unlikely that
the endomorphism $\sigma^p = \Hom(\hat{\psi}^p, K_\wp^\times)$ of 
$K_1(K_\wp \Gammam)$ induces an endomorphism of $K_0T(\okpg)$ ($p:=
\char(\cO_K/\wp)$) in general.

\bigskip

\section*{\S 4 The Galois Structure of the Module of Relative \\
Differentials}

Let $N/K$ be a Galois extension of number fields with Galois group
$\Gammam = \Gal(N/K)$. Let $\Omega := \Omega^1_{\cO_N/\cO_K}$ denote the
module of relative differentials, $\wD := {\rm Ann}_{\cO_N}(\Omega)$ the 
different, $V:= {\rm supp}(\Omega)$ the set of ramified primes in 
$\Spec(\cO_N)$, and $\wJ   $ the $\Gammam$-stable ideal $\prod_{\wP \in V} \wP$
in $\cO_N$. 

The aim of this section is to prove the following theorem which will be
used in section 5 to compute the tangential element $T_f$ associated with
the morphism $f: \Spec(\cO_N) \ra \Spec(\cO_K)$.

{\bf Theorem 4.1.}  We have:
\[[\Omega] = [\wJ   ] - [\wJ   \wD] \quad {\rm in} \quad K'_0(\cO_N \fis  \Gammam).\]

{\bf Proof.} The quotient module $\wJ  /\wJ  \wD$ decomposes into the direct
sum of the $\Gammam$-modules $\wJ  _\wp/\wJ  _\wp \wD_\wp$, $\wp \in f(V)$, over $\cO_N$.
Thus, we have:
\[[\wJ  ]-[\wJ  \wD] = \sum_{\wp \in f(V)} [\wJ  _\wp/\wJ  _\wp \wD_\wp] \quad {\rm in}
\quad K'_0(\cO_N \fis  \Gammam).\]
Likewise, we have:
\[[\Omega] = \sum_{\wp \in f(V)} [\Omega_\wp] \quad {\rm in} \quad
K'_0(\cO_N\fis  \Gammam).\]
Therefore, it suffices to show that, for all $\wp \in f(V)$, we have
\[[\Omega_\wp] = [\wJ  _\wp/\wJ  _\wp \wD_\wp]\]
in the Grothendieck group $K_0T(\Gammam, \cO_{N,\wp})$ of all f.~g.\ 
torsion modules over $\cO_{N,\wp}\fis \Gammam$. \\
Now we fix a $\wp \in f(V)$ and a $\wP \in V$ over $\wp$. Let $\Gammam_\wP
\subseteq \Gammam$ denote the decomposition group of $\wP$. Since the
canonical homomorphism $K_0T(\cO_{N,\wp}\fis \Gammam) \ra K_0T(\cO_{N,\wP}\fis 
\Gammam_\wP)$ is bijective, it suffices to show the equality
\[[\Omega_\wP] = [\wP\cO_{N,\wP}/\wP\wD_\wP] \quad {\rm in} \quad 
K_0T(\Gammam_\wP, \cO_{N,\wP}).\]
We now write $\Gammam$, $\cO_N$, $\wP$, $\Omega$, and $\wD$ for
$\Gammam_\wP$, $\cO_{N,\wP}$, $\wP\cO_{N,\wP}$, $\Omega_{\wP}$, and
$\wD_\wP$, respectively. Let $\pi \in \wP$ be a prime element and let
$L\subseteq N$ denote the inertia field associated with $\wP$. We again write
$\cO_L$ for the localization of $\cO_L$ in $\wP \cap L$. By 
Proposition 18 on p.\ 19 in \cite{Se2}, the $\cO_L$-algebra $\cO_N$
is generated by $\pi$. Thus, the module $\Omega_{\cO_N/\cO_L}$ of
relative differentials is generated by $d\pi$. Since $\cO_L$ is 
unramified over $\cO_{K,\wp}$, the canonical homomorphism 
$\Omega_{\cO_N/\cO_L} \ra \Omega$ is bijective. Hence, $\Omega$ is
generated by $d\pi$, too. Let the unit $u_\gamma \in \cO_N^\times$
be defined by $\gamma(\pi) = u_\gamma \cdot \pi$ for each $\gamma \in 
\Gammam$. Then, for all $i \ge 1$, we have $\gamma(\pi^i) = u_\gamma^i 
\cdot \pi^i$ and, hence,
\[\gamma(\pi^i \, d\pi) = u_\gamma^{i+1} \pi^{i} \, d\pi +
u_\gamma^i \pi^{i+1} \, du_\gamma \quad {\rm in} \quad \Omega.\]
Consequently, we have:
\[\gamma(\overline{\pi^i \, d\pi}) = u_\gamma^{i+1} 
\overline{\pi^i \, d\pi} \quad {\rm in} \quad \wP^i\Omega/\wP^{i+1}\Omega.\]
Thus, the sequence
\[\begin{array}{ccccccccc}
0& \ra & \wP^{i+2} & \stackrel{\rm can}{\ra} & \wP^{i+1} & \ra & 
\wP^i \Omega/\wP^{i+1}\Omega & \ra & 0\\
&&&& a \cdot \pi^{i+1} & \mapsto & a \cdot \overline {\pi^i \, d\pi}
\end{array}\]
($a \in \cO_N$) is an exact sequence of $\Gammam$-modules over $\cO_N$ for all
$i = 0, \ldots, l-1$ (here, $l$ denotes the length of $\Omega$, i.e.,
$\wD = \wP^l$). This implies the desired equality
\[[\Omega] = \sum_{i=0}^{l-1} [\wP^i\Omega/\wP^{i+1}\Omega] =
\sum_{i=1}^{l} [\wP^i/\wP^{i+1}] = [\wP/\wP\wD] \quad {\rm in} \quad
K_0T(\Gammam, \cO_N).\]

{\bf Corollary 4.2.} If $N/K$ is a tame extension, we have
\[[\Omega] = [\wD^{-1}] - [\cO_N] \quad {\rm in} \quad K'_0(\cO_N\fis \Gammam).\]

{\bf Proof.} By Theorem (2.6) on p.\ 210 in \cite{N}, we have
\[\wJ  \wD = \prod_{\wP \in V} \wP^{e_\wP}\]
where $e_\wP$ denotes the ramification index at $\wP$. Thus, for all
$\wp \in f(V)$, the $\Gammam$-stable ideal $\wJ  _\wp \wD_\wp$ is equal to
$\wp \cO_{N, \wp}$ and then isomorphic to $\cO_{N,\wp}$ as a $\Gammam$-module
over $\cO_{N,\wp}$. Thus, we have:
\[[\Omega] = \sum_{\wp \in f(V)} [\Omega_\wp] = \sum_{\wp \in f(V)}
[\Omega_\wp \otimes (\wJ  \wD)_\wp^{-1}] = [\Omega \otimes (\wJ  \wD)^{-1}]
\quad {\rm in} \quad K_0T(\Gammam, \cO_N).\]
By Theorem 4.1, we finally have:
\[[\Omega] = [\Omega \otimes (\wJ  \wD)^{-1}] = ([\wJ  ]-[\wJ  \wD]) \cdot ([\wJ  \wD]^{-1})
= [\wD^{-1}] - [\cO_N] \quad {\rm in} \quad K'_0(\cO_N\fis \Gammam).\]

The following Remark 4.3 and Example 4.4 contain situations where
the equality $[\Omega] = [\wD^{-1}] - [\cO_N]$ holds not only in the
Grothendieck group $K'_0(\cO_N\fis \Gammam)$, but $\Omega$ and $\wD^{-1}/\cO_N$
are in fact equivariantly isomorphic.

{\bf Remark 4.3.} Let $\wP \in \Spec(\cO_N)$ be a tamely ramified prime
ideal and let $I\subseteq \Gammam$ denote the inertia group associated with
$\wP$. Then the modules $\Omega_\wP$ and $\wD_\wP^{-1}/\cO_{N,\wP}$
are isomorphic as $I$-modules over $\cO_{N,\wP}$. 

{\bf Proof.} Let $L:=N^I$ denote the inertia field associated with $\wP$,
$\wp := \wP \cap L$ the prime below $\wP$, and $e$ the ramification index
at $\wP$. Furthermore, let $\hat{\cO}_{N,\wP}$ and $\hat{\cO}_{L,\wp}$
denote the respective completions. 
Then there is a prime element $\pi \in \hat{\cO}_{N,\wP}$
such that $\pi^e \in \hat{\cO}_{L,\wp}$ (e.g., see p.\ 26 in \cite{F1}).
As in the proof of Theorem 4.1, the module of relative differentials
$\Omega_\wP \cong \Omega_{\cO_{N,\wP}/\cO_{L,\wp}} \cong 
\Omega_{\hat{\cO}_{N,\wP}/\hat{\cO}_{L,\wp}}$ is generated by $d\pi$. 
Let the $e$-th root of unity $\zeta_\gamma$ defined by $\gamma(\pi) =
\zeta_\gamma \cdot \pi$ for each $\gamma \in I$. Since
$\zeta_\gamma$ is contained in $\hat{\cO}_{L,\wp}$ (e.g., see
p.\ 26 in \cite{F1}), we have:
\[\gamma(d\pi) = \zeta_\gamma \, d\pi + \pi \, d\zeta_\gamma = \zeta_\gamma
\, d\pi \quad {\rm in} \quad \Omega_\wP.\]
Hence, the sequence
\[\begin{array}{ccccccccc}
0 & \ra & \wP^e \hat{\cO}_{N,\wP} & \ra & \wP \hat{\cO}_{N,\wP} & \ra 
& \Omega_\wP & \ra & 0\\
&&&& a \cdot \pi & \mapsto & a \cdot d\pi
\end{array}\]
($a \in \cO_N$) is an exact sequence of $I$-modules over $\hat{\cO}_{N,\wP}$. The 
homomorphism $\wD_\wP^{-1}\ra \wP\hat{\cO}_{N,\wP}$, $a \mapsto
a \cdot \pi^e$, of $I$-modules over $\cO_{N,\wP}$ now induces the desired
isomorphism
\[\wD_\wP^{-1}/\cO_{N,\wP} \,\, \tilde{\ra} \,\, \Omega_\wP\]
of $I$-modules over $\cO_{N,\wP}$. 

{\bf Example 4.4.} Let $D$ be a square-free natural number, $K$ the
field of rational numbers, and $N:= \QQ(\sqrt{D})$. Then the modules $\Omega$
and $\wD^{-1}/\cO_N$ are isomorphic as $\Gammam$-modules over $\cO_N$. 

{\bf Proof.} We have
\[\cO_N = \left\{ \begin{array}{rll}
\ZZ[\sqrt{D}]  \cong  \ZZ[T]/(\phi) &{\rm where} \; \phi = T^2 - D, &
{\rm if} \; D\not\equiv 1 \; {\rm mod} \; 4\\
\ZZ[(1+\sqrt{D})/2]  \cong  \ZZ[T]/(\phi) & {\rm where} \; 
\phi = T^2 - T - (D-1)/4,& {\rm if} \; D \equiv 1 \; {\rm mod} \; 4.
\end{array} \right.\]
Let $t:= \sqrt{D}$ or $t:= (1+\sqrt{D})/2$. By Satz (2.4) on p.\ 207 in
\cite{N}, the different $\wD_{N/K}$ is generated by $\phi'(t)$ which
is equal to $2\sqrt{D}$ or $\sqrt{D}$. Let $\gamma \in \Gammam =\Gal(N/K)$ denote
the non-trivial automorphism. Then we have $\gamma(t) = -t$ or 
$\gamma(t) = -t +1$ and hence
\[\gamma(dt) = -dt \quad {\rm in} \quad \Omega\]
in both cases. Consequently, the sequence
\[\begin{array}{ccccccccc}
0 & \ra & \cO_N & \stackrel{\rm can}{\ra} & \wD^{-1} & \ra & 
\Omega & \ra & 0\\
&&&& a \cdot (\phi'(t)^{-1}) & \mapsto & a \cdot dt
\end{array}\]
($a \in \cO_N$) is an exact sequence of $\Gammam$-modules over $\cO_N$. This proves 
Example 4.4.

\bigskip

\section*{\S 5 Equivariant Adams-Riemann-Roch Formulas for \\
Tame Galois Extensions}

Let $N/K$ be a {\em tame} Galois extension of number fields with 
Galois group $\Gammam = \Gal(N/K)$. Let $f: \Spec(\cO_N) \ra \Spec(\cO_K)$
denote the $\Gammam$-morphism associated with $N/K$. As in section 3, let
$\psi_k^{\rm CNT}$ denote the endomorphism of the locally free 
classgroup ${\rm Cl}(\okg)$ defined by Ph.\ Cassou-Nogu\`es and M.\
Taylor in \cite{CNT}.

Any $\Gammam$-stable ideal $\wA$ in $N$ defines an element $(\wA) \in 
{\rm Cl}(\cO_K)$ in
the obvious way (see below). In the paper \cite{BC}, D.\ Burns and
T.\ Chinburg have established a formula for $\psi_k^{\rm CNT}((\wA))$. 
In this section, we will give a geometric interpretation of this 
formula: On the one hand, we will apply the equivariant Adams-Riemann-Roch
formula of \cite{KoGRR} to the $\Gammam$-morphism $f$ using Corollary 4.2.
On the other hand, we will reformulate the formula of Burns and Chinburg
using Theorem 3.7. It will then turn out, that, more or less, the
formula of Burns and Chinburg is a strengthening of the equivariant
Adams-Riemann-Roch formula in this situation.

Since $\cO_N$ and $\cO_K$ are regular rings, the morphism $f$ is
a local complete intersection morphism (see Remark on p.\ 86 in \cite{FL}).
Being a finite morphism, $f$ in particular is a projective morphism.
By Remark (3.5) in \cite{KoGRR}, the morphism $f$ then is 
$\Gammam$-projective in the sense of section 3 in \cite{KoGRR}. Let 
$T^\vee_f \in K_0(\Gammam, \cO_N)$ denote the equivariant cotangential
element associated with $f$. It is defined as follows (see section 4
in \cite{KoGRR}): Let 
\[\Spec(\cO_N) \,\, \stackrel{i}{\longrightarrow} \,\,
\AA = \Spec(A) \,\, \stackrel{p}{\longrightarrow} \,\,
\Spec(\cO_K)\]
be a factorization of $f$ into a closed $\Gammam$-immersion $i$ and
a smooth $\Gammam$-morphism $p$ of constant relative dimension $d$. 
Then $i$ is a regular embedding (see Remark on p.\ 86 in \cite{FL}).
Let $I$ denote the $\Gammam$-stable ideal in $A$ associated with $i$.
Then $T^\vee_f$ is defined to be
\[T_f^\vee := [ \Omega^1_{A/\cO_K} \otimes_A A/I] - [I/I^2] \in
K_0(\Gammam, \cO_N).\]

{\bf Lemma 5.1.} We have:
\[T_f^\vee = [\Omega^1_{\cO_N/\cO_K}] \quad {\rm in} \quad K'_0(\cO_N\fis 
\Gammam).\]

{\bf Proof.} We have a natural exact sequence
\[I/I^2 \ra \Omega^1_{A/\cO_K} \otimes_A A/I \,\,
\stackrel{\epsilon}{\ra} \,\, \Omega^1_{\cO_N/\cO_K} \ra 0 \]
of $\Gammam$-modules over $\cO_N$ (e.g., see Theorem 58i) on p.\ 187 
in \cite{Mat}).
Being a submodule of $\Omega^1_{A/\cO_K} \otimes A/I$, the kernel of
$\epsilon$ is torsionfree; thus, it is $\cO_N$-projective. Furthermore,
we obviously have $\rank_{\cO_N}(\ker(\epsilon)) = d = \rank_{\cO_N}(I/I^2)$.
Thus, the canonical surjection $I/I^2 \ra \ker(\epsilon)$ is an 
isomorphism. This implies that the left homomorphism in the sequence above
is injective which proves Lemma 5.1.

Let $\wD$ denote the different associated with $N/K$.

{\bf Corollary 5.2.} We have:
\[T_f^\vee = [\wD^{-1}] - [\cO_N] \quad {\rm in} \quad K_0(\Gammam, \cO_N).\]

{\bf Proof.} This immediately follows from Lemma 5.1 and Corollary 4.2
since $K_0(\Gammam, \cO_N) \cong K'_0(\cO_N\fis \Gammam)$ (see Notations).

Now, we fix a $k\in \NN$. Let $\hat{K}_0(\Gammam, \cO_K)[k^{-1}]$ denote
the $J$-adic completion of $K_0(\Gammam, \cO_K)[k^{-1}]$ where
$J:= \ker(\rank: K_0(\Gammam, \cO_K) \ra \ZZ)$ is the augmentation ideal.

{\bf Proposition 5.3.} The Bott element $\theta^k(\wD^{-1}) := 
\sum_{i=0}^{k-1}[\wD^{-1}]$ is invertible in 
\[K_0(\Gammam, \cO_N) 
\otimes_{K_0(\Gammam, \cO_K)} \hat{K}_0(\Gammam,\cO_K)[k^{-1}]\] 
and we have:
\[\theta^k(T^\vee_f)^{-1} = k \cdot \theta^k(\wD^{-1})^{-1}\]
(see section 4 in \cite{KoGRR} for the definition of $\theta^k(T^\vee_f)^{-1}$).

{\bf Proof.} For any f.~g.\ $\cO_N$-projective $\cO_N\fis \Gammam$-module $P$, 
let $\theta^k(P) \in K_0(\Gammam,\cO_N)$ denote the
$k$-th Bott element associated with $P$ (see section 4 in \cite{KoGRR}). 
By section
4 in \cite{KoGRR}, there is a representation $T^\vee_f = [\Omega] - [\cN]$
in $K_0(\Gammam, \cO_N)$ such that $\theta^k(\Omega)$ is invertible in
$K_0(\Gammam, \cO_N) \otimes_{K_0(\Gammam, \cO_K)} \hat{K}_0(\Gammam, \cO_K)
[k^{-1}]$. By Corollary 5.2, we have $\theta^k(\wD^{-1}) \cdot 
\theta^k(\cN) = k \cdot \theta^k(\Omega)$. Thus, $\theta^k(\wD^{-1})$,
$\theta^k(\cN)$, and $\theta^k(T^\vee_f)$ are invertible and we have
$\theta^k(T^\vee_f)^{-1} = k\cdot \theta^k(\wD^{-1})^{-1}$. 

Pulling back the module structure along $\cO_K \ra \cO_N$ defines a
homomorphism 
\[f_*: K_0(\Gammam, \cO_N) \ra K_0(\Gammam, \cO_K).\]
Obviously, $f_*$ coincides with the push-forward homomorphism $f_*$
defined in section 3 of \cite{KoGRR}. The projection formula implies
that $f_*$ induces a homomorphism
\[\hat{f}_* : K_0(\Gammam,\cO_N) \otimes_{K_0(\Gammam, \cO_K)} 
\hat{K}_0(\Gammam, \cO_K)[k^{-1}] \ra \hat{K}_0(\Gammam, \cO_K)[k^{-1}].\]
Let $\psi^k$ denote the $k$-th Adams operation on $K_0(\Gammam, \cO_K)$ and
$K_0(\Gammam, \cO_N)$.

{\bf Theorem 5.4} (Equivariant Adams-Riemann-Roch formula applied to $f$).
For all $x\in K_0(\Gammam, \cO_N)$, we have:
\[\psi^k f_* (x) = \hat{f}_*(k\cdot \theta^k(\wD^{-1})^{-1} \cdot \psi^k(x))
\quad {\rm in} \quad \hat{K}_0(\Gammam, \cO_K)[k^{-1}].\]

{\bf Proof.} Because of Proposition 5.3, this is the assertion of
Theorem (4.5) in \cite{KoGRR} applied to $f$.

We are now going to reformulate the formula of Burns and Chinburg
mentioned in the beginning of this section. 

{\bf Lemma 5.5.} \\
(a) Every f.\ g.\ $\cO_N$-projective $\cO_N\fis\Gammam$-module has a
quotient which is isomorphic to a $\Gammam$-stable fractional ideal in $N$.
In particular, any invertible $\cO_N\fis\Gammam$-module is isomorphic to
such an ideal.\\
(b) Every f.\ g.\ $\cO_N$-projective $\cO_N\fis\Gammam$-module is 
projective as $\cO_K\Gammam$-module.\\
(c) The (additive) group $K_0(\Gammam, \cO_N)$ is generated by the
classes of $\Gammam$-stable fractional ideals in $N$. 

{\bf Proof.} Obviously, any f.\ g.\ $\cO_N$-projective 
$\cO_N\fis\Gammam$-module $P$ is a $\Gammam$-stable lattice in 
$P\otimes _{\cO_N} N$. By Morita theory, the $N\fis \Gammam$-module 
$P \otimes_{\cO_N} N$ is isomorphic to  $N^{\rank(P)}$. Thus, there
is an $N\fis\Gammam$-submodule $W$ of $P\otimes_{\cO_N} N$ such that
$(P \otimes_{\cO_N} N)/W$ is isomorphic to $N$. Then, $P/P\cap W$ is
isomorphic to a fractional ideal in $N$. This shows 
assertion (a). Using (a), one easily shows as in \S 3 of chapter I in
\cite{F1} that $P$ is locally free over $\cO_K\Gammam$. By Theorem 1.1 on
p.\ 1 in \cite{T}, this proves assertion (b). Assertion (c) follows 
from assertion (a) by induction on the rank of $P$. 

By Lemma 5.5(b), we have a push-forward homomorphism $f_*: K_0(\Gammam,\cO_N)
\ra K_0(\cO_K\Gammam)$ such that the following diagram commutes:
\[\begin{array}{ccccc}
&& K_0(\Gammam, \cO_N)\\
\\
& f_* \swarrow && \searrow f_*\\
\\
K_0(\cO_K\Gammam) && \stackrel{c}{\longrightarrow}&& K_0(\Gammam, \cO_K).
\end{array}\]
Let ${\rm Ind}_1^\Gammam ({\rm Cl}(\cO_K))$
denote the image of the classical classgroup ${\rm Cl}(\cO_K) \subset 
K_0(\cO_K)$ under the induction map
${\rm Ind}_1^\Gammam: K_0(\cO_K) \ra K_0(\cO_K\Gammam)$. 

We now assume that the prime divisors of $\ord(\Gammam)$ are unramified
in $K$ and that $k$ is coprime to the exponent $e(\Gammam)$ of $\Gammam$. 
We fix a $k' \in \NN$
such that $kk' \equiv 1$ mod $e(\Gammam)$.

{\bf Theorem 5.6} (Reformulation of a formula of Burns and Chinburg).
We have for all $x\in K_0(\Gammam,\cO_N)$:
\[\psi^k f_*(x) = f_* \left(k \cdot \sum_{i=0}^{k'-1} [\wD^{-ik}] \cdot
\psi^k(x)\right) \quad {\rm in} \quad K_0(\cO_K \Gammam)/
\left({\rm Ind}_1^\Gammam({\rm Cl}(\cO_K))\oplus e(\Gammam)\ZZ[\cO_K\Gammam]
\right).\]
If the classical classgroup ${\rm Cl}(\cO_K)$ is trivial and if
the rank of $x$ is $0$, then this formula already holds in $K_0(\cO_K\Gammam)$.

{\bf Proof.} Let $\psi_{k'}^{\rm CNT}$ denote the endomorphism of
${\rm Cl}(\cO_K\Gammam)$ defined by Cassou-Nogu\`es and Taylor (e.g., see
section 3). For any $\Gammam$-stable fractional ideal $\wA$ in $N$,
$(\wA) := [f_*(\wA)] - [\cO_K\Gammam]$ is a well-defined element of
${\rm Cl}(\cO_K \Gammam)$ by Lemma 5.5(b). By Corollary 2.7 in \cite{BC}, 
we have:
\[\psi^{\rm CNT}_{k'} ((\wA)) = \sum_{i=0}^{k'-1} (\wD^{-ik} \cdot
\wA^k) \quad {\rm in} \quad {\rm Cl}(\cO_K\Gammam)/{\rm Ind}_1^\Gammam 
({\rm Cl}(\cO_K)).\]
Using Theorem 3.7, we obtain the equality
\[\psi^k([f_*(\wA)]-[\cO_K\Gammam]) = k \cdot \sum_{i=0}^{k'-1} 
([f_*(\wD^{-ik}\cdot \wA^k)] - [\cO_K\Gammam]) \quad {\rm in} \quad
K_0(\cO_K\Gammam)/{\rm Ind}_1^\Gammam({\rm Cl}(\cO_K)).\]
Hence, we get the following equality
in $K_0(\cO_K\Gammam)/{\rm Ind}_1^\Gammam({\rm Cl}(\cO_K))$
for all $x \in K_0(\Gammam, \cO_N)$ (by Lemma 5.5(c)):  
\[\psi^k f_*(x) = f_*\left(k\cdot \sum_{i=0}^{k'-1} 
[\wD^{-ik}]\cdot \psi^k(x)\right)
- kk' \cdot \rank_{\cO_N}(x) \cdot [\cO_K\Gammam] + \rank_{\cO_N}(x) \cdot
\psi^k([\cO_K\Gammam]).\]
Now, Theorem 1.6(e) shows Theorem 5.6.

Using an idea of the proof of Theorem 2 in \S 4 in \cite{CEPT}, we show
in the following lemma that the right side of the formula in 
Theorem 5.4 coincides with the right side of the formula in 
Theorem 5.6 (without $f_*$) for elements of the form $x= T_f^\vee \cdot y$ with
$y \in K_0(\Gammam, \cO_N)$.

{\bf Lemma 5.7.} We assume that the $\cO_N\fis\Gammam$-modules
$\wD^{\ord(\Gammam)}$ and $\cO_N$ are isomorphic. (This certainly holds if
the classical classgroup ${\rm Cl}(\cO_K)$
is trivial.) Then we have for all $y \in K_0(\Gammam,\cO_N)$:
\[\theta^k(\wD^{-1})^{-1} \cdot \psi^k(T_f^\vee \cdot y) =
\sum_{i=0}^{k'-1} [\wD^{-ik}] \cdot \psi^k(T_f^\vee \cdot y) \quad
{\rm in} \quad K_0(\Gammam, \cO_N) \otimes_{K_0(\Gammam, \cO_K)}
\hat{K}_0(\Gammam, \cO_K)[k^{-1}].\]

{\bf Proof.} We have:
\begin{eqnarray*}
\lefteqn{\theta^k(\wD^{-1})^{-1} \cdot \psi^k(T_f^\vee \cdot y) =}\\
&& = (\sum_{i=0}^{k-1} [\wD^{-i}])^{-1} \cdot 
([\wD^{-k}] - [\cO_N]) \cdot \psi^k(y) \quad 
\mbox{(by Corollary 5.2)}\\
&& = ([\wD^{-1}] - [\cO_N]) \cdot \psi^k(y) \quad
\mbox{(geometric series)}\\
&& = ([\wD^{-kk'}] - [\cO_N])\cdot \psi^k(y) \quad \mbox{
(since $\wD^{\ord(\Gammam)} \cong \cO_N$ by assumption)}\\
&& = (\sum_{i=0}^{k'-1} [\wD^{-ik}]) \cdot ([\wD^{-k}]-[\cO_N]) \cdot
\psi^k(y) \quad \mbox{(geometric series)}\\
&& = \sum_{i=0}^{k'-1} [\wD^{-ik}] \cdot \psi^k(T_f^\vee \cdot y)
\quad \mbox{(by Corollary 5.2)}.
\end{eqnarray*}

{\bf Corollary 5.8.} We assume that ${\rm Cl}(\cO_K)$ is trivial. Then
the element $f_*(T_f^\vee)$ in $K_0(\cO_K\Gammam)$ or $K_0(\Gammam, \cO_K)$
is an eigenvector of the Adams operation $\psi^k$ with eigenvalue $k$.

{\bf Proof.} This follows from Theorem 5.6 and the proof of Lemma 5.7.

In the following closing remark, we compare Theorems 5.4 and 5.6 and
point out their respective advantages.

{\bf Remark 5.9.} On the one hand, Theorem 5.4 is a special case of a
general (geometric) equivariant Adams-Riemann-Roch formula which 
can be applied also in the non-tame case, which does not need the
assumption ``$\gcd(k,\ord(\Gammam)) =1$'', and which holds without passing
to residues modulo ${\rm Ind}_1^\Gammam(K_0(\cO_K))$. On the other hand,
Theorem 5.6 holds already for $K_0(\cO_K\Gammam)$ without passing to
$\hat{K}_0(\Gammam, \cO_K)[k^{-1}]$ or even $K_0(\Gammam, \cO_K)$. 
In particular, the
multiplier $\sum_{i=0}^{k'-1} [\wD^{-ik}]$ in Theorem 5.6 is definable
already in $K_0(\Gammam, \cO_N)$. In fact, Theorem 5.6 is not a 
reformulation but a weakening of the formula of Burns and Chinburg 
since we have multiplied their formula by the factor $k$ in the proof
of Theorem 5.6. Finally, choosing $k=1$ and $k'= 1+e(\Gammam)$ in Theorem
5.6 yields the (non-trivial!) formula 
\[f_*\left(\sum_{i=1}^{e(\Gammam)} [\wD^{-i}]\cdot x\right) =0 \quad
{\rm in} \quad K_0(\cO_K\Gammam)\left/\left({\rm Ind}_1^{\Gammam}({\rm Cl}(\cO_K))
\oplus e(\Gammam) \ZZ[\cO_K\Gammam]\right)\right.\]
for all $x\in K_0(\Gammam, \cO_N)$ (see also Remark 2.8 in \cite{BC}). One 
can show that, after passing to $\hat{K}_0(\Gammam, \cO_K)[k^{-1}]
\left/\left({\rm Ind}_1^{\Gammam}({\rm Cl}(\cO_K)) \oplus e(\Gammam)
\ZZ[\cO_K\Gammam]\right)\right.$, this formula together with the formula
in Theorem 5.4 is in fact equivalent to the general form of the formula
in Theorem 5.6.\\ 
Moreover, the paper \cite{BC} contains
a formula for $\psi_{k}^{\rm CNT}((\wA))$ even if $\gcd(k,\ord(\Gammam))
\not= 1$ (see Corollary 2.4 in \cite{BC}). Since, however, there is
no relation of $\psi_k^{\rm CNT}$ with Adams operations in this case
(see Remarks 3.11 and 3.12), this formula cannot be compared with 
the Adams-Riemann-Roch formula.

\bigskip

Mathematisches Institut II\\
Universit\"at Karlsruhe\\
D-76128 Karlsruhe\\
Germany

{\em e-mail:} Bernhard.Koeck@math.uni-karlsruhe.de


\begin{thebibliography}{CEPT}
\bibitem[AB]{AB} {\sc  K.\ Akin} and {\sc  D.\ A.\ Buchsbaum}, 
Characteristic-free representation theory of the general linear group, 
{\em Adv.\ Math.}\  {\bf 58} (1985), 149-200.
\bibitem[ABW]{ABW} {\sc  D.\ Akin, D.\ A.\ Buchsbaum,} and {\sc  J.\ Weyman}, 
Schur functors and Schur complexes, {\em Adv.\ Math.}\  {\bf 44} 
(1982) 207-278.
\bibitem[B]{B} {\sc  H.\ Bass}, ``Algebraic $K$-theory'', {\em Math.\  
Lecture Note Series} (Benjamin, New York, 1968).
\bibitem[BC]{BC} {\sc  D.\ Burns} and {\sc  T.\ Chinburg}, Adams operations 
and integral Hermitian-Galois representations, {\em Amer.\ J.\ Math.}\  
{\bf 118} (1996), 925-962.
\bibitem[CNT]{CNT} {\sc  Ph.\ Cassou-Nogu\`es} and {\sc  M.\ J.\ Taylor}, 
Op\'erations d'Adams et Groupe des classes d'Alg\`ebre de groupe, 
{\em J.\ Algebra} {\bf 95} (1985), 125-152.
\bibitem[C]{C} {\sc  T.\ Chinburg}, Galois structure of de Rham cohomology
of tame covers of schemes, {\em Ann.\ of Math.}\  {\bf 139} (1994), 443-490.
\bibitem[CEPT]{CEPT} {\sc  T.\ Chinburg}, {\sc  B.\ Erez}, {\sc  G.\ 
Pappas,} and {\sc  M.\ J.\ Taylor}, Arithmetic equivariant Riemann-Roch theorems,
{\em preprint} (Manchester Centre for Pure Mathematics, number 1993/16).
\bibitem[CR]{CR} {\sc  C.\ W.\ Curtis} and {\sc  I.\ Reiner}, ``Methods of 
representation theory with applications to finite groups and orders'',
vol.\ I, {\em Pure Appl.\ Math.}\ (Wiley, New York, 1981).
\bibitem[F1]{F1} {\sc  A.\ Fr\"ohlich}, ``Galois module structure of algebraic 
integers'', {\em Ergeb.\ Math.\ Grenz\-geb.\ (3)}
{\bf 1} (Springer, New York, 1983).
\bibitem[F2]{F2} {\sc  A.\ Fr\"ohlich}, ``Classgroups and Hermitian modules'',
{\em Progr.\ Math.}\  {\bf 48} (Birk\-h\"au\-ser, Boston, 1984).
\bibitem[FH]{FH} {\sc W.\ Fulton} and {\sc J.\ Harris}, ``Representation
theory: a first course'', {\em Grad.\ Texts in Math.}\ {\bf 129} (Springer,
New York, 1991).
\bibitem[FL]{FL} {\sc  W.\ Fulton} and {\sc  S.\ Lang}, 
``Riemann-Roch algebra'', 
{\em Grundlehren Math.\ Wiss.}\ {\bf 277} (Springer, New York, 1985).
\bibitem[Gr1]{Gr1} {\sc  D.\ R.\ Grayson}, Higher algebraic $K$-theory: II,
in {\sc  R.\ M.\ Stein} (ed.), ``Algebraic $K$-theory (Evanston, 1976)'',
{\em Lecture Notes in Math.}\ {\bf 551} (Springer, New York, 1976), 217-240.
\bibitem[Gr2]{Gr2} {\sc  D.\ R.\ Grayson}, Exterior power operations on higher
$K$-theory, {\em $K$-Theory} {\bf 3} (1989), 247-260.
\bibitem[Ke]{Ke} {\sc  M.\ Kervaire}, Op\'erations d'Adams en Th\'eorie des
repr\'esentations lin\'eaires des groupes finis, {\em Enseign.\ Math.}\ 
{\bf 22} (1976), 1-28.
\bibitem[Ko1]{KoARR} {\sc  B.\ K\"ock}, Das Adams-Riemann-Roch-Theorem in der 
h\"oheren \"aquivarianten $K$-Theorie, {\em J.\ Reine Angew.\ Math.}\
{\bf 421} (1991), 189-217.
\bibitem[Ko2]{KoSh} {\sc  B.\ K\"ock}, Shuffle products in higher $K$-theory,
{\em J.\ Pure Appl.\ Algebra} {\bf 92} (1994), 269-307. 
\bibitem[Ko3]{KoAdHi} {\sc  B.\ K\"ock}, On Adams operations 
on the higher $K$-theory of group rings, in:
G.\ Banaszak et al.\ (eds.), ``Algebraic $K$-theory (Pozna\'n, 1995)'', 
{\em Contemp.\ Math.}\ {\bf 199} (Amer.\ Math.\ Soc., Providence, 1996), 
139-150.
\bibitem[Ko4]{KoAdPro} {\sc  B.\ K\"ock}, Adams operations for projective
modules over group rings, {\em Math.\ Proc.\ Cambridge 
Philos.\ Soc.} {\bf 122} (1997), 55-71.
\bibitem[Ko5]{KoGRR} {\sc  B.\ K\"ock}, The Grothendieck-Riemann-Roch theorem
for group scheme actions, {\em Ann.\ Sci.\ \'Ecole Norm.\ Sup.}\ {\bf 31}
(1998), 415-458.
\bibitem[McD]{MacD} {\sc  I.\ G.\ Macdonald}, ``Symmetric functions and Hall
polynomials'', {\em Oxford Math.\ Monographs} (Clarendon Press, Oxford, 1979).
\bibitem[Ma]{Mat} {\sc  H.\ Matsumura}, ``Commutative algebra'', {\em Math.\ 
Lecture Note Series} (Benjamin, New York, 1970).  
\bibitem[N]{N} {\sc  J.\ Neukirch}, ``Algebraische Zahlentheorie'' (Springer,
Berlin, 1992).
\bibitem[Que]{Que} {\sc  J.\ Queyrut}, $S$-groupes des classes d'un ordre
arithm\'etique, {\em J.\ Algebra} {\bf 76} (1982), 234-260.
\bibitem[Q]{Q} {\sc  D.\ Quillen}, Higher algebraic $K$-theory: I, in 
{\sc  H.\ Bass} (ed.), ``Algebraic $K$-Theory I (Seattle, 1972)'', 
{\em Lecture Notes in Math.}\ {\bf 341} (Springer, New York, 1973), 85-147.
\bibitem[Se1]{Se1} {\sc  J.\ P.\ Serre}, ``Lineare Darstellung endlicher
Gruppen'', {\em Logik und Grundlagen der Mathematik} (Vieweg, Braunschweig, 1972).
\bibitem[Se2]{Se2} {\sc  J.\ P.\ Serre}, ``Local fields'', {\em Grad.\ Texts in 
Math.}\  {\bf 67} (Springer, New York, 1979).
\bibitem[Sh]{Sh} {\sc C.\ Sherman}, Connecting homomorphisms in localization
sequences, in: 
G.\ Banaszak et al.\ (eds.), ``Algebraic $K$-theory (Pozna\'n, 1995)'', 
{\em Contemp.\ Math.}\ {\bf 199} (Amer.\ Math.\ Soc., Providence, 1996), 
175-183.
\bibitem[Sw]{Sw} {\sc  R.\ G.\ Swan}, Induced representations and projective
modules, {\em Ann.\ of Math.}\  {\bf 71} (1960), 552-578.
\bibitem[T]{T} {\sc  M.\ Taylor}, ``Classgroups of group rings'', {\em Lecture
Notes Ser.}\ {\bf 91} (Cambridge University Press, Cambridge, 1984).

\end{thebibliography}
\end{document}